\documentclass{amsart}[article]
\usepackage{amssymb}
\usepackage{amsmath}
\usepackage[foot]{amsaddr}
\usepackage{bbm}
\usepackage{bbold}
\usepackage{bm}
\usepackage{braket}
\usepackage[title]{appendix}
\usepackage[dvipsnames]{xcolor}
\usepackage{soul}
\usepackage{graphicx}
\usepackage{textcomp}
\usepackage{enumerate}
\usepackage{verbatim}
\usepackage{enumitem}
\usepackage{mathrsfs}
\usepackage[normalem]{ulem}
\usepackage{caption}

\usepackage{hyperref}
\hypersetup{
    colorlinks,
     urlcolor={blue!80!black}
}
\usepackage{lipsum, babel}
\graphicspath{{../figura1/}{../figura2/}}
\usepackage{graphicx}

\usepackage{bbold}
\usepackage{bm}
\usepackage{a4wide}

\usepackage{soul}

\newcommand{\RE}{\mathbb{R}}

\newcommand{\BX}{\mathbb{X}}
\newcommand{\YY}{\mathbb{Y}}
\newcommand{\PP}{\mathbb{P}}
\newcommand{\II}{\mathbb{I}}

\newcommand{\UU}{\mathbb U}
\newcommand{\SSS}{\mathfrak S}
\newcommand{\HH}{\mathfrak H}
\newcommand{\VV}{\mathbb{V}}

\newcommand{\bx}{\mathbf{x}}
\newcommand{\by}{\mathbf{y}}
\newcommand{\ee}{\mathbf{e}}
\newcommand{\bs}{\mathbf{s}}
\newcommand{\CS}{\mathcal{S}}
\newcommand{\CF}{\mathcal{F}}

\newcommand{\Cov}{\mathsf{Cov}}

\newcommand{\CC}{\mathcal{C}}
\newcommand{\CX}{\mathcal{X}}
\newcommand{\CA}{\mathcal{A}}

\definecolor{Teal1}{rgb}{0.0, 0.5, 0.5}
\definecolor{DarkGreen1}{rgb}{0.0, 0.28, 0.15}
\newcommand{\la}{\lambda}

\newcommand{\bz}{\mathbf{z}}
\newcommand{\bzero}{\mathbf{0}}

\newcommand{\CN}{\mathcal{N}}

\usepackage[scr]{rsfso}
\newcommand{\CK}{\mathcal{K}}
\newcommand{\CKK}{\mathscr{K}}
\newcommand{\oK}{{K}}
\newcommand{\bK}{{\mathsf{K}}}

\newcommand{\s}{\sigma}

\newcommand{\E}{\mathbf{E}}

\DeclareMathOperator{\tr}{tr}
\DeclareMathOperator{\vvec}{vec}

\usepackage{xcite}
\usepackage{xr}

\newcommand{\bS}{\mathbf{S}}

\newcommand{\CL}{\mathcal{L}}
\newcommand{\CLL}{\CL^{+,s}_1}

\newcommand{\CBB}{\mathscr{B}}

\newcommand{\bY}{\mathbf{Y}}

\newcommand{\eI}{\mathbb1}

\newtheorem{definition}{Definition}[section]
\newtheorem{theorem}[definition]{Theorem}
\newtheorem{lemma}[definition]{Lemma}

\newtheorem{proposition}[definition]{Proposition}

\newtheorem{remark}{Remark}

\newtheorem{example}{Esempio}

\def\CC{\mathcal C}
\def\E{\mathbb E}

\def\P{\mathbb P}

\def\been{\begin{enumerate}}
\def\bee{\begin{example}}
\def\bel{\begin{lemma}}
\def\bepr{\begin{proposition}}
\def\bep{\begin{proof}}
\def\bet{\begin{theorem}}

\def\enen{\end{enumerate}}
\def\ene{\end{example}}
\def\enl{\end{lemma}}
\def\enpr{\end{proposition}}
\def\enp{\end{proof}}
\def\ent{\end{theorem}}

\def\la{\lambda}

\def\s{\sigma}

\DeclareMathAlphabet{\mathsfit}{T1}{\sfdefault}{\mddefault}{\sldefault}

\title{
LDP for the  covariance process  in fully connected neural networks.
}

\author{Luisa Andreis}
\author{Federico Bassetti}
\author{Christian Hirsch}
\address[Luisa Andreis]{Department of Mathematics, Politecnico di Milano, Milan, Italy.}
\address[Federico Bassetti]{Department of Mathematics, Politecnico di Milano, Milan, Italy.}
\address[Christian Hirsch]{Department of Mathematics\\ Aarhus University \\ Ny Munkegade, 118, 8000, Aarhus, Denmark.}
\address[Christian Hirsch]{DIGIT Center, Aarhus University, Finlandsgade 22, 8200 Aarhus, Denmark}
\email{luisa.andreis@polimi.it}
\email{federico.bassetti@polimi.it}
\email{hirsch@math.au.dk}
\date{\today}  % Inserisce automaticamente la data di oggi

\begin{document}
\begin{abstract}
	  In this work, we study large deviation properties of the covariance process in fully connected Gaussian deep neural networks. More precisely, we establish a large deviation principle (LDP) for the covariance process in a functional framework, viewing it as a process in the space of continuous functions. As key 	applications of our main results, we obtain  posterior LDPs under Gaussian likelihood in both the infinite-width and mean-field regimes. The proof is based on an LDP for the covariance process as a Markov process valued in the space of non-negative, symmetric trace-class operators equipped with the trace norm.
\end{abstract}

\maketitle


\section{Introduction}
In large neural networks with many neurons per layer, key theoretical insights emerge through the study of scaling limits, particularly the infinite-width limit, wherein the network depth remains fixed while the number of neurons per layer tends to infinity. In this regime, the behavior of neural networks simplifies significantly, and Gaussian processes arise as the limiting behavior.

This Gaussian universality appears both in training under gradient flow, as described by the neural tangent kernel \cite{JacotNTK}, and in the Bayesian inference setting, where precise correspondences between neural networks and kernel methods have been established \cite{lee2018deep, g.2018gaussian}. 
In the Bayesian setting, 
this phenomenon was first established in Neal’s seminal work \cite{Neal} for shallow networks and has since been extended to fully connected networks with nonlinear activations \cite{lee2018deep, g.2018gaussian, Hanin2023}, as well as to certain convolutional architectures \cite{novak2019bayesian, garriga-alonso2018deep}.

Gaussian processes arises naturally due to central limit effects in the network’s outputs. When weights are chosen to be gaussian, a key observation  is that  the output of the network (at finite size) is a mixture of Gaussians with a random covariance, which, as width increases, converges to a deterministic covariance function. The output thus converges in distribution to a Gaussian process,  known as the Neural Network Gaussian Process (NNGP).

Once the deterministic  limit of the covariance is identified, a natural next step is to study fluctuations around this limit, particularly through the lens of large deviations. Recent works \cite{MaPaTo24, Vog24} have developed large and moderate deviation principles for the rescaled output of fully connected Gaussian networks with a finite input set. In point of fact, in  these studies, the large deviation principle (LDP) is first established for the random covariance matrix, and then transferred to the network's output.

The present paper extends this perspective by formulating, in a Bayesian framework, both a law of large numbers and a large deviation principle at the functional level for the random covariance function under Gaussian prior.  Once the large deviation principle (LDP) for the covariance under the prior is established, a corresponding result under the posterior—assuming a Gaussian likelihood—follows relatively straightforwardly.

In this setting, the posterior distribution reflects the training of the neural network. In contrast, in the non-Bayesian setting, training is typically performed via optimization using stochastic gradient descent. The latter introduces intricate dependencies, making large deviation analysis considerably more challenging and requiring entirely different mathematical tools. Initial progress in this direction has been made for networks with a single hidden layer \cite{whnn}.

Turning to the methods employed, our approach centers on modeling the vector of covariances as a Markov process in the natural space of non-negative, symmetric trace-class operators. This functional-analytic perspective allows us to derive our main results with conceptual clarity. In particular, we invoke LDP  for sums of independent (but non-identically distributed) Banach space-valued random variables \cite[Theorem~2]{Bolthausen1984}.
Here, a crucial ingredient in establishing exponential tightness is the Gaussianity of the weights. Indeed, we combine the classical  Fernique theorem with a useful compactness criteria developed by \cite{deAcosta85}. Finally, we extend our results to the sup-norm topology, which is particularly relevant in neural network applications. This is again achieved via exponential tightness arguments, leveraging on a concentration result for the largest eigenvalue of Gaussian matrices, borrowed from  \cite{Ver19}.

To conclude this introduction, we highlight an interesting  implication of our findings. The simplification of random neural networks in the NNGP regime comes at a notable cost: a significant loss in the model’s expressive and learning capabilities. This limitation is particularly evident when contrasted with modern deep architectures, which are capable of rich feature learning well beyond the capabilities of networks in the infinite-width regime \cite{ChizatLazy, lewkowycz2021the, NEURIPS2020_1457c0d6}.
Here, we show that the LDP under the posterior coincides with that of the prior, leading to an identical rate function. This phenomenon reflects the laziness of the infinite-width regime: the training data does not affect the rate function, underscoring the absence of genuine learning in this limit.

However, this Gaussian—and effectively “lazy”—behavior is not unavoidable. Alternative regimes, yielding different dynamics, have been proposed in literature:
 the mean field  scaling, see e.g.  \cite{doi:10.1073/pnas.1806579115,https://doi.org/10.1002/cpa.22074,doi:10.1137/18M1192184,NEURIPS2018_a1afc58c},
the heavy tailed initial weight distributions setting, see e.g.  \cite{BordinoFavaroFortini2023,FavaroFortiniPeluchetti2023} or 
the so-called proportional limit (where both the number of training patterns $P$ and the number of neurons $N$ diverges at the same rate) investigated
in physics literature   \cite{pacelli2023statistical,aiudi2023,baglioni2024predictive}.  

In our work, following \cite{BassettietalJMLR}, we focus on the mean-field parametrization and prove that, in this regime, the posterior distribution of the covariance function satisfies a LDP with a non-trivial rate function, that explicitly reflects the influence of the training data.

To summarize, the main 	contributions of this paper are the following:

\begin{enumerate}
    \item In \textbf{Theorem~\ref{thm.main_C0}} we prove an LDP for the covariance process in a functional setting, that is as a process in the space of continuous functions.
    \item The key step is \textbf{Theorem~\ref{thm.main}}, where we prove the LDP for the covariance process as a Markov process taking values
    in the space of \emph{non-negative} and \emph{symmetric trace-class  operators} with trace norm.
    \item A key application area for our main results are \textbf{ Proposition~\ref{thm.main_C0_posteriod}} and \textbf{Proposition~\ref{thm.main_C0_posteriodMF}},
    where LDP under the posterior - with Gaussian likelihood -- is derived 
    both for the usual infinite-width limit and for the mean--field parametrization. 
\end{enumerate}

The rest of the manuscript is organized as follows. In Section \ref{sec:NN} we introduce the considered neural network model together with the associated covariance process. Then, we state our main functional LDP results in Section \ref{sec:ldp}. The LDP on the posterior process of Bayesian neural networks  is presented in Section \ref{sec:post}. The main building blocks for our proofs are outlined in Section~\ref{sec:main_ingredients}, the main steps of the proofs are in Section \ref{Proofs_mainsteps}, whereas all the missing details are given in Section \ref{sec:proofs}. Finally, we include an appendix discussing basic properties of the considered function space topology and containing some of the highly technical proof steps. In particular, in Appendix~\ref{appendix:notation}, to ease the reading, Table~\ref{tab:notations} collects the important pieces of notation.

\section{The Neural network and the covariance process}
\label{sec:NN}

In this section we introduce the main definitions and review some well-known facts 
on Neural Networks. 

\subsection{Fully-connected Bayesian  deep  neural networks}

In a  {\it fully-connected neural network with $L$ hidden layers}, 
the pre-activations of each layer $h^{(\ell)}=(h^{(\ell)}_1,\dots,h^{(\ell)}_{N_\ell})$   are given recursively as functions of the pre-activations of the previous layer $h^{(\ell-1)}=(h^{(\ell-1)}_1,\dots,h^{(\ell-1)}_{N_{\ell-1}})$:
for a given input $\mathbf x=(x_1,\dots,x_{N_0})$
in $\RE^{N_0}$
\begin{equation}
    \label{main_recursion}
    \begin{split}
h_i^{(1)}(\mathbf x) &= \frac1{\sqrt {N_{0}}} \sum_{j=1}^{N_{0}} W^{(0)}_{i j} x_j \, 
\quad i=1,\dots,N_1,
\\
h_i^{(\ell)}(\mathbf x) &= \frac1{\sqrt {N_{\ell-1}}} \sum_{j=1}^{N_{\ell-1}} W^{(\ell-1)}_{i j} \s\Big ( h_j^{(\ell-1)}(\mathbf x) \Big ) \,
\quad i=1,\dots,N_\ell, \quad \ell\ge 2,   \\
\end{split}
\end{equation}
where $W^{(\ell-1)}$ are  the weights and we assume that the so-called  biases of the $\ell$-th layer are zero. 
Assuming the output has $N_{L+1}=D$ components,  the function implemented by the neural network   is  
 the output of the last layer $h^{(L+1)}(\bx|\vartheta  )=(h_1^{(L+1)}(\mathbf x),\dots,h_{D}^{(L+1)}(\mathbf x))^\top$.
  Here 
\[
\vartheta = \{W^{(\ell-1)}_{ij}: \ell=1,\dots,L+1;  i=1,\dots, N_\ell, j=1,\dots, N_{\ell-1} \}
\] 
represents the collection of all the trainable weights of the network.

In a Bayesian neural network, 
 a prior for the weights $\vartheta $ is specified, which 
translates in a prior for $h^{(\ell)}(\bx)=h^{(\ell)}(\bx|\vartheta  )$. 
%%%% 
Hence, given a compact subset $\UU \subset \RE^{N_0}$ of inputs,  
 the corresponding output at layer $\ell$ is 
 the stochastic vector field  
\[
[h^{(\ell)}(\mathbf x)]_{\mathbf x \in \UU}=
[(h_1^{(\ell)}(\mathbf x), \dots,h_{N_\ell}^{(\ell)}(\mathbf x))^\top]_{\mathbf x \in \UU}. 
\]
In the following we shall consider the  lazy-training  infinite  width  limit, where the depth is fixed, and the width (number of neurons) grows to infinity, and 
a gaussian prior for the  for the weights.

In summary, we shall assume the following hypotheses:
\begin{enumerate}[label=(H\arabic*)]
\item \label{H1} \textbf{Gaussian prior}: 
 the  weights are Gaussian, i.e. 
 \begin{equation*}\label{lawofW}
   W^{(\ell)}_{ij}  \stackrel{ind}{\sim}  \CN\big(0,\la_\ell^{-1}\big)
\quad    \ell=1,\dots,L+1;  \, i=1,\dots, N_\ell, \,j=1,\dots, N_{\ell-1} ,
 \end{equation*}
 with $\la_\ell^{-1}\in (0,\infty)$ layer dependent variance.
\item \label{H2}
\textbf{Polynomial growth condition}:  the activation function $\s: \RE \to \RE$ is continuous and there is  $r>0$
and $A<+\infty$
such that 
 \[
 \s(x)^2 \le A (1+| x|^{r}).
 \]  
\item \label{H3} \textbf{Infinite width regime}: $N_0$ and $N_{L+1}=D$ are fixed, while  
$N_\ell=N_\ell(N)$ with   $N \to +\infty$
and $\lim_{N\to +\infty} N_\ell/N=m_\ell \in (0,+\infty)$.
\end{enumerate} 
 Here, and in the rest of the paper $\CN(\mathbf{m},\mathbf{C})$ denotes the Gaussian distribution 
 with mean $\mathbf{m}$ and covariance matrix $\mathbf{C}$. 
We assume that the random weights are defined on a common probability space $(\Omega,\CF,\PP)$.

For some results we shall need also  
the following additional hypothesis.

\begin{enumerate}[label=(H\arabic*),start=4] 
    \item \label{H4} \textbf{Lipschitz condition}: the function $\s$ is uniformly Lipschitz continuous, i.e. there is  $L_{\s}<\infty$ such that 
\[
(\s(x)-\s(y))^2\le L^2_{\s}(x-y)^2,
\]
for all $x,y$.
\end{enumerate}

\begin{remark}\label{rem:bias1} 
   Various works consider as well an independent bias at each step of the recursion. 
 To ease notation and exposition, we choose not to do so (or equivalently we choose this bias to be identically equal to zero at every layer). All our results hold in presence of Gaussian bias as well, we give more details in Remark~\ref{rem:bias2}.
   \end{remark}
\subsection{The Neural Network Gaussian Process}\label{conditional_covariance_finDIM}

Under assumptions \ref{H1}-\ref{H2},  conditionally on the penultimate layer, the output of the network is Gaussian. In other words, the output (under the prior) is a mixture of Gaussian processes with random covariance functions. Such mixtures are generally difficult to characterize, but in the infinite-width limit, the random covariance converges to a deterministic limit. Consequently, the output converges in distribution to a pure Gaussian process. This limiting behavior is commonly referred to as the Neural Network Gaussian Process (NNGP) limit, which we now briefly review.

Let
    \[
  \CF^{\ell}=\s\Big (h_i^{(\ell)}(\bx): \bx \in \UU, i=1,\dots,N_{\ell} \Big )
  \] 
the $\s$-field containing all the information up to layer $\ell$,  where  $\CF^0$ is the trivial $\s$-field. 
Simple computations show that:

\noindent\textit{conditionally on   $\CF^{\ell}$, for any finite input set $\UU_P=\{\bx_1,\dots,\bx_P\}$, the collection of random variables 
$[ h^{(\ell)}_i(\bx)]_{ i=1,\dots,N_\ell,\, \, \bx \in \UU_P}$ are jointly normal  with zero mean and  conditional covariances 
\begin{equation}\label{mainCovariance}
\Cov( h^{(\ell)}_i(\bx_\mu), h^{(\ell)}_j(\bx_\nu) | \CF^{\ell-1})=
  \delta_{ij} \CK^{\ell}_{N_{\ell-1}}(\bx_\mu,\bx_\nu),
\end{equation}
for $i,j=1,\dots, N_{\ell}$ and $\mu,\nu=1,\dots, P$, with  
  \begin{equation}
      \label{eq:K_in_C0}  \begin{split}
   \CK^1_{N_0}(\bx,\bx') & \colon =\frac1{ \la_0  {N_{0}}  }  \sum_{i_{0}=1}^{N_0 }   \mathbf x_{i_{0}} \mathbf x_{i_{0}}'   \\
 \CK^{\ell}_{N_{\ell-1}}(\bx,\bx')& \colon = \frac1{   \la_{\ell-1} {N_{\ell-1}}} \sum_{i=1}^{N_{\ell-1}} 
\s\Big ( h_i^{(\ell-1)}(\mathbf x) \Big )\s\Big ( h_i^{(\ell-1)}(\mathbf x') \Big ) \quad \text{for $\ell=2,\dots,L+1$}.
\end{split}
  \end{equation}}
It is easy to see that  $\omega \mapsto h^{(\ell)}(\mathbf x)$ is measurable and that $\mathbf x \mapsto h^{(\ell)}(\mathbf x)$ 
belongs to $C^0(\UU,\RE^{N_\ell})$, the space 
of continuous function on $\UU$ with values in $\RE^{N_\ell}$. 
 Hence  $h^{(\ell)}$ are jointly measurable and 
$\omega \mapsto  \{ h_i^{(\ell)}(\mathbf x): \bx \in \UU\}$ is a random element with values in $C^0(\UU,\RE)$ for every $i$.
By 
 well-known properties of Gaussian processes
(see, e.g., Lemma 13.1 \cite{KallenbergBook}), 
using \eqref{mainCovariance}, 
one proves the following important fact. 

\begin{enumerate}[label=(F\arabic*)]
\item \label{FactGP}
\textit{Assuming \ref{H1}-\ref{H2}, conditionally on   $\CF^{\ell-1}$,  the stochastic processes  $\{h^{(\ell)}_i(\bx): \bx \in \UU\}$ {\rm(}$i=1,\dots,N_\ell${\rm)}  are independent and Gaussian 
with zero mean and  conditional covariance function $\CK^{\ell}_{N_{\ell-1}}$
given in \eqref{eq:K_in_C0}.}
\end{enumerate}

This fact was first pointed out in \cite{novak2019bayesian}. At this stage the key observation is that,  given a finite input set  $\UU_P=\{\bx_1,\dots,\bx_P\}$, 
the  sequence of random covariance matrices 
\[
C^{\ell,P}_{N_{\ell-1}}\colon =[ \CK^{\ell}_{N_{\ell-1}}(\bx,\bx')]_{(\bx,\bx') \in \UU_P^2}
\]
is a Markov chain with deterministic initial condition $C^{1,P}_{N_{0}}$  and it 
converges  to a sequence of  deterministic matrices, satisfying the NNGP recursion. 

The NNGP recursion sets   the initial condition
 $\CK^1_{\infty}(\bx,\bx')\colon =\CK^1_{N_0}(\bx,\bx')$  and defines 
for $\ell=2,\dots,L+1$ the covariance functions 
  \[
  \CK^{\ell}_{\infty}(\bx,\bx') \colon = \frac1{   \la_{\ell-1}}  \E \Big [ 
\s\Big ( h^{\ell-1}_\infty(\mathbf x) \Big )\s\Big ( h^{\ell-1}_{\infty} (\mathbf x') \Big )  \Big ] \quad  \text{with} \quad
 h^{\ell-1}_{\infty} \sim GP(\mathbf{0},\CK^{\ell-1}_\infty) 
\]
where $Z \sim GP(\mathbf{0},\CKK)$ 
means that $Z$ is a Gaussian process with zero mean and covariance $\CKK$, i.e. 
$\CKK(\bx,\bx')=\E[Z(\bx)Z(\bx')]$. 
Note that, since $\s$ is polynomially bounded, the $ \CK^{\ell}_{\infty}(\bx,\bx')$'s are well-defined continuous covariances for any $\ell$.

\textbf{Covariance concentration:} 
\textit{assuming \ref{H1}-\ref{H2}-\ref{H3}, 
 given 
$\UU_P=\{\bx_1,\dots,\bx_P\}$, one has
\begin{equation}\label{LLN_COV}
(C^{2,P}_{N_1}\dots,C^{L+1,P}_{N_L}) \stackrel{\PP}{\to} ({C}^{2,P},\dots,{C}^{L+1,P}) \quad \text{as $N \to \infty$}
\qquad \text{\rm(law of large numbers)}
\end{equation}
where
$
{C}^{\ell,P}=[ \CK^{\ell}_{\infty}(\bx_\mu,\bx_\nu)]
_{(\bx_\nu,\bx_\mu) \in \UU_P^2}$.
}

The previous result is contained, more or less explicitly, in various works, see in particular 
Lemma 2.3 in \cite{Hanin2023} and
Proposition 21 in \cite{hanin2024}. 
The random covariance is an instance of what is  called a {\it collective observable} in \cite{hanin2024}.

Due to the fact that the random covariance matrix converges to a deterministic limit, 
the neural networks
simplify significantly in the infinite width regime, as summarized by the next important result, where
the random processes $h^{(\ell)}_i$ implicitly depend on $N$.

\textbf{Gaussian limit:}
\textit{assuming \ref{H1}-\ref{H2}-\ref{H3}, 
given 
$\UU_P=\{\bx_1,\dots,\bx_P\}$, 
at each layer and given any integer $M<+\infty$, when  $N \to +\infty$,
\begin{equation}\label{CLT_P}
[h^{(\ell)}_i ( \bx) ]_{\bx \in \UU_P;i=1,\dots,M} \stackrel{\CL}{\to} \CN(0,\mathbf{C}^{\ell})
\qquad \text{\rm(central limit)}
\end{equation}
where
$
\mathbf{C}^{\ell}_{{(i,\mu),(j,\nu)}} \colon =\delta_{ij} [{C}^{\ell,P}]_{\mu,\nu}
$.
}

The result above has  been obtained many times and under a variety of different assumptions,
including more general network architectures. See, e.g. \cite{lee2018deep,g.2018gaussian,Hanin2023}. 
We refer the interested reader to \cite{Hanin2023}
 for a discussion. 
A functional central limit theorem is proved in 
\cite{Hanin2023}.  
Convergence rates to Gaussian limit for fully connected networks  have been derived  in \cite{favaro2023quantitative, Trevisan2023}. In~\cite{MaPaTo24, Vog24} large and moderate deviations for the output vectors $[h^{(L+1)} ( \bx_\mu) ]_{\mu=1,\dots,P}$, properly rescaled by $1/\sqrt{N}$, are obtained. 
Various results for the very special case of 
deep linear networks, e.g. can be found in \cite{doi:10.1073/pnas.2301345120,zavatone-veth2021exact,Mufan2023,BaLaRo}, in particular we mention  that 
\cite{BassettietalJMLR} proves  a LDP for the covariance structure of  deep linear fully connected networks. 

See Section~\ref{sec:literature} for a deeper comparison between some of these works and our results.

\section{LDP for the covariance process}
\label{sec:ldp}
In this Section, we state our main results, namely the functional LDPs for the covariance process, see Theorem \ref{thm.main} and Theorem~\ref{thm.main_C0} below. To introduce the suitable function space, we  provide the necessary functional analytic preliminaries in Section \ref{ss:trace_class}. Next, to motivate our LDP, we first state the LLN  in Section \ref{ss:LLN}. Section \ref{Main.thm} contains the statement of the LDPs, respectively on the space of trace-class operators and on the space of continuous functions. They hold under slightly stronger assumptions than the LLN, as it often happens. In Section~\ref{sec:literature} we perform a short overview on related literature, comparing our results to previous ones in similar frameworks.

\subsection{Trace-class operators and Gaussian random elements in Hilbert spaces} 
\label{ss:trace_class}
Given     a separable  Hilbert space $H$  with scalar product $(\cdot,\cdot)_H$, 
we denote  by  $\CL_1(H)$ the Banach space  of  trace-class operators   on $H$, endowed with the trace  norm    
$\|\oK\|_1\colon =\tr|\oK|$, where  $|\oK|=\sqrt{\oK^*\oK}$. 
The closed cone of \emph{non-negative} and \emph{symmetric} trace-class  operators will be denoted by  $\CL_1^{+,s}(H)$. See Appendix~\ref{app:trace_class} for details.

A  measure $\gamma$ on $\CBB(H)$ (the Borel $\s$-field on $H$) is said to be  a
\textit{Gaussian  of zero mean 
and covariance $\oK \in \CL_1^{+,s}(H)$}, if 
$\gamma \circ g^{-1} =\CN(0,\sigma^2)$ with $\sigma^2=(\oK g,g)_H$ for every $g \in H$.
In what follows we denote by $\CN_H(\bzero,\oK)$  such a measure. 
In particular, for a random element $Z$ taking values in $(H,\CBB(H))$, 
 we write  $Z \sim \CN_H(\bzero,\oK)$, 
if for every $g \in H$ the real-valued  random variable $(g,Z)_H$ is  Gaussian with zero mean and variance $(\oK g,g)_H$.  See  \cite{Bogachev} for further details.

It will be important the following fact:
\begin{enumerate}[label=(F\arabic*),start=2]
\item \label{Fact1}
 {\it Given a sequence of operators
$\oK_n \in \CL_1^{+,s}(H)$ and an operator $\oK \in \CL_1^{+,s}(H)$ one has 
that $ \CN_H(\bzero,\oK_n)$ converges weakly as a measure to $ \CN_H(\bzero,\oK)$ if and only if 
$\|\oK_n -\oK\|_1 \to 0$. }  
\end{enumerate}
To see this combine    Ex. (iii) 3.8.13  in  \cite{Bogachev} 
with Lemma \ref{Lemma_eqvivalent_conv}
in  Appendix. 

In what follows $L^2(\UU)$ denotes the space of (measurable)  functions $f: \UU \to \RE$ 
such that $\|f\|_{L^2}^2 \colon = \int_\UU |f(x)|^2 dx<+\infty$.

If  $\{Z(\bx,\omega)\}_{\bx\in\UU}$ is a (jointly measurable) Gaussian process defined on a compact set $\UU \subset \RE^{N_0}$ with paths in $H=L^2(\UU)$, zero mean and covariance function $\CKK$, 
then it  can be seen as Gaussian random element with values in $H$ and   $Z \sim \CN_H(\bzero,\oK)$
where   $\oK$ 
   is the covariance operator canonically  associated to $\CKK$ through 
\begin{equation}\label{mapKerneltoCov}
\oK g(\bx) = \int_{\UU}  \CKK (\bx,\by) g(\by) d\by, \qquad g \in H.
\end{equation}
See Examples 2.3.16 and  3.11.14  in \cite{Bogachev}.

In the rest of the manuscript 
we fix
\[
H\colon =L^2(\UU), \quad 
 (f,g)_H=\int_\UU f(\bx)g(\bx) d\bx, \quad
\CL_1=\CL_1(L^2(\UU)) \quad  \text{and} \quad  
\CL_1^{+, s}=\CL_1^{+, s}(L^2(\UU)).
\]
Finally, 
let $\CC^{+,s} \subset C^0(\UU^2,\RE)$ be the  class of continuous, symmetric, 
positive definite kernels 
on $\UU^2$ (see \ref{Mercer} in Appendix). 
By \eqref{mapKerneltoCov} we define 
\[
\phi(\CKK)\colon =\oK \quad 
\phi: \CC^{+,s} \to \CL_1^{+,s}.
\]
This is a well-defined and continuous map, 
 see next Lemma \ref{continuity_Phi}.
 
\subsection{The law of large numbers for the covariance process in the space of trace-class operators} 
\label{ss:LLN}
Since $\UU$ is compact, $C^0(\UU,\RE) \subset L^2(\UU)$ with continuity. It follows 
that $\omega \mapsto  \{ h^{(\ell)}_i(\mathbf x): \bx \in \UU\}$   can be regarded as a  random element with values in $H=L^2(\UU)$  
equipped with its Borel $\s$-field. 
Analogous considerations hold for 
$ \omega \mapsto \{ \CK^{\ell}_{N_{\ell-1}}(\bx,\bx'):  (\bx,\bx') \in \UU^2\}  $, where now one needs to consider $C^0(\UU \times \UU,\RE)$
and $L^2(\UU \times \UU)$. 
Note also that, being $\CK^\ell_{N_{\ell-1}}$ continuous covariance kernels, 
$\P(\CK^\ell_{N_{\ell-1}} \in \CC^{+,s})=1$.

Now, introduce 
 the (random) trace-class operators on $H$
defined by 
\[
[\bK^\ell_{N_{\ell-1}} g](\bx) \colon =(\CK^{\ell}_{N_{\ell-1}}(\bx,\cdot),g)_H 
\]for all $g\in H$, i.e.\ $\bK^\ell_{N_{\ell-1}} =\phi(\CK^{\ell}_{N_{\ell-1}})$. For $\ell=1$ the operator $\bK^1_{N_0}=\phi(\CK_{N_0}^1)$ is deterministic  
and does not depend on $N$, while 
for $\ell= 2,\dots,L+1$ one has 
\begin{equation}\label{conditional_covarianceOp}
[\bK^\ell_{N_{\ell-1}} g](\bx) =(\CK^{\ell}_{N_{\ell-1}}(\bx,\cdot),g)_H =\frac1{   \la_{\ell-1} {N_{\ell-1}}} \sum_{i=1}^{N_{\ell-1}}
\s\Big ( h_i^{(\ell-1)}(\mathbf x) \Big )
 \int_{\UU}
\s\Big ( h_i^{(\ell-1)}(\mathbf y) \Big )  g(\by) d\by.
\end{equation}
To be more formal, for any $f \in H$, introduce  the linear operator on $H$ defined by
\begin{equation}\label{def_Cf}
[C_fg](\bx)= \s(f(\bx)) \int_{\UU} \s(f(\by)) g(\by)d\by \qquad \forall g \in H. 
\end{equation}
 Note that
 $f \mapsto C_f$ is a continuous application
 from $H$ to $\CL_1^{+,s}$
 (by the next Lemma \ref{Lemma_cont_Cf}). 

With these definitions,  after  identifying  $h_i^{(\ell)}$ with random elements taking  values in the Hilbert space $H$, thanks to 
\ref{FactGP},
we can state that for $\ell \ge 1$ 
\begin{equation}\label{conditional_h}
h_i^{(\ell)}  \big | \CF^{\ell-1} \stackrel{iid}{\sim } \CN_{H}(\bzero,\bK^{\ell}_{N_{\ell-1}})  \quad  i=1,\dots,N_{\ell}.
\end{equation}

As in the finite dimensional case, it is also easy to see that the sequence of  covariance random operators  $\bK^2_{N_1},\dots,\bK^{L+1}_{N_L}$ is a Markov chain with deterministic initial condition $\bK^1_{N_0}$. 
Combining \eqref{conditional_covarianceOp} and
\eqref{conditional_h}
the transition rule  can be  described by 
\begin{equation}\label{transitionK}
\bK^{\ell}_{N_{\ell-1}}  \mapsto  \bK^{\ell+1}_{N_\ell} \colon =
\frac1{   \la_\ell N_\ell}  \sum_{i=1}^{N_\ell} C_{ h_i^{(\ell)}}
   \quad
      h_i^{(\ell)}  \big | \bK^{\ell}_{N_{\ell-1}} \stackrel{iid}{\sim } \CN_{H}(\bzero,\bK^{\ell}_{N_{\ell-1}})  \quad  i=1,\dots,N_\ell . \\
\end{equation}
See Lemma \ref{regualr_CondProb_kernel} for an explicit  construction of the corresponding  transition kernel. 

At this stage, note that the NNGP recursion
induces a sequence of  
 covariance operators $\bK_\infty^\ell\colon =\phi(\CK^{\ell}_{\infty})$.   
In analogy with \eqref{LLN_COV}, we have the following result. 

\begin{proposition}[LLN in $\CLL$] \label{prop:LLN_traceclass}
Assume \ref{H1}-\ref{H2}-\ref{H3},
then 
\begin{equation}\label{LLN_K}
(\bK^2_{N_0},\dots,\bK^{L+1}_{N_L}) \stackrel{\PP}{\to}  (\bK^2_{\infty},\dots,\bK^{L+1}_{\infty}) \quad \text{as $N \to \infty$}
\qquad \text{(law of large numbers)}.
\end{equation}
\end{proposition}

We are now ready to discuss the large deviations 
of 
 $(\bK^2_{N_1},\dots,\bK^{L+1}_{N_L})$
 and 
$(\CK^2_{N_1},\dots, \CK^{L+1}_{N_L})$.

\subsection{Large deviation of the covariance process in $\CLL$ and in $\CC^{+,s}$}\label{Main.thm}

Denote by  $\CL_\infty$ the set of bounded linear operators on $H$. It is well-known that the dual of  $\CL_1$ can be  isometrically identified to $\CL_\infty$, with duality product given by  $\langle D, C \rangle =\tr(DC)$ for $D \in \CL_\infty$ and $C \in \CL_1$. See \ref{item:sep_cond} in the Appendix. 

Fix $\lambda>0$ and define, for any $K_1$ and $K_2$  in $\CL_1^{+,s}$, the following function 
\begin{equation}\label{eq:rate_function}
I_\la(K_2|K_1)\colon =\sup_{ D \in \CL_\infty} \Big \{ \tr(DK_2) 
- \log\Big( \int_H e^{ \frac{1}{\la}\tr(DC_h)}
 \CN_H(dh|\bzero,K_1) \Big )\Big \}.
\end{equation}

\begin{theorem}[LDP in $\CLL$]\label{thm.main}
Assume that $\UU \subset \RE^{N_0}$ is compact and  \ref{H1}-\ref{H2}-\ref{H3}, with
$r<2$ in \ref{H2}.
Then,  the sequence of laws of $\{(\bK^2_{N_1},\dots,\bK^{L+1}_{N_L})\}_N$
 satisfies     the LDP on 
        $\CLL \times \dots \times \CLL$ with speed $N$ and good rate 
function 
\[
I(K_2,\dots,K_{L+1})\colon ={m_1}I_{\la_1}(K_2|\bK^1_{N_0})
+ {m_2}I_{\la_2}(K_3|K_2)+\cdots+{m_L}I_{\la_L}(K_{L+1}|K_L),
\]
for all $(K_2,\dots,K_{L+1})$ in $\CLL \times \dots \times \CLL$. 
\end{theorem}

The above LDP holds for slightly stronger assumptions than the LLN in Proposition~\ref{LLN_K}, indeed we require $r<2$ in \ref{H2}. We expect the same result to hold for $r=2$ as well, but proving this would require to deal with significant technical details (see~\cite{Vog24} for the finite dimensional case) and we postpone it to future work.

From the above LDP,  requiring some additional regularity of $\sigma$ (assumption~\ref{H4}), one can as well obtain a LDP on the space $\CC^{+,s}\times \dots\times \CC^{+,s}$.

\begin{theorem}[LDP in $\CC^{+,s}$]\label{thm.main_C0}
Assume that $\UU \subset \RE^{N_0}$ is compact and  \ref{H1}-\ref{H2}-\ref{H3}-\ref{H4}, with
$r<2$ in \ref{H2}.
Then, the sequence of laws of $\{(\CK^2_{N_1},\dots, \CK^{L+1}_{N_L})\}_N$
 satisfies     the LDP on 
    $\CC^{+,s} \times \dots \times \CC^{+,s}$ with speed $N$ and good rate 
function 
\[
\mathcal{I}(\CKK_2,\dots,\CKK_{L+1})\colon=I(\phi(\CKK_2),\dots,\phi(\CKK_{L+1}))
\]
for all $(\CKK_2,\dots,\CKK_{L+1})\in \CC^{+,s} \times \dots \times \CC^{+,s}$. 
\end{theorem}

Our approach is built on the idea of representing the vector of covariances as a Markov process whose values lie in a particularly natural mathematical setting: the space of non-negative, symmetric, trace-class operators.
This choice of state space is not only conceptually fitting, but also analytically powerful—it allows us to leverage a rich body of probabilistic tools. In particular, many of the limiting results we seek follow directly from limit theorems for sums of independent (though not necessarily identically distributed) random variables in Banach spaces,~\cite[Theorem~2]{Bolthausen1984}.
To handle the dependence introduced by the Markovian structure, we apply a conditional version of these limit theorems, which, together with exponential tightness, enables us to iteratively construct the overall limit by combining the results obtained at each step of the chain. The full proofs are  in Section~\ref{Proofs_mainsteps}.

\begin{remark}\label{rem:bias2}
    We mention in Remark~\ref{rem:bias1} that often this model includes a bias, that is the addition of an independent Gaussian random variable to the recursion. Here we give more details. Adding the bias consists in substituting \eqref{main_recursion} with the following
   \begin{equation*}
    \begin{split}
h_i^{(1)}(\mathbf x) &= B^{(1)}_i+ \frac1{\sqrt {N_{0}}} \sum_{j=1}^{N_{0}} W^{(0)}_{i j} x_j \, 
\quad i=1,\dots,N_1
\\
h_i^{(\ell)}(\mathbf x) &=B^{(\ell)}_i+ \frac1{\sqrt {N_{\ell-1}}} \sum_{j=1}^{N_{\ell-1}} W^{(\ell-1)}_{i j} \s\Big ( h_j^{(\ell-1)}(\mathbf x) \Big ) \,
\quad i=1,\dots,N_\ell, \quad \ell\ge 2,   \\
\end{split}
\end{equation*}
where for any $\ell=1,\dots,L+1$, $\mathbf{B}^{(\ell)}=(B^{(\ell)}_i)_{i=1,\dots,N_{\ell}}$ is a $N_\ell$-dimensional vector of i.i.d. Gaussian with zero mean and variance $b^{(\ell)}\geq 0$, independent from the weights, i.e.
\begin{equation*}
   B^{(\ell)}_{i}  \stackrel{ind}{\sim}  \CN\big(0,b^{(\ell)}\big),\quad \ell=1,\dots, L+1;\, i=1,\dots, N_{\ell}.
 \end{equation*}
       This results in a simple deterministic translation of the covariance function defined in~\eqref{eq:K_in_C0}, that becomes
     \begin{equation*}
     \begin{split}
   \CK^1_{N_0}(\bx,\bx') & \colon =b^{(1)}+\frac1{ \la_0  {N_{0}}  }  \sum_{i_{0}=1}^{N_0 }   \mathbf x_{i_{0}} \mathbf x_{i_{0}}'   \\
 \CK^{\ell}_{N_{\ell-1}}(\bx,\bx')& \colon = b^{(\ell)}+ \frac1{   \la_{\ell-1} {N_{\ell-1}}} \sum_{i=1}^{N_{\ell-1}} 
\s\Big ( h_i^{(\ell-1)}(\mathbf x) \Big )\s\Big ( h_i^{(\ell-1)}(\mathbf x') \Big ) \quad \text{for $\ell=2,\dots,L+1$}.
\end{split}
  \end{equation*}
  Consequently \eqref{conditional_covarianceOp} is substituted by
\begin{equation}\label{conditional_covarianceOp_with_bias}
[\bK^\ell_{N_{\ell-1}} g](\bx) =b^{(\ell)}(1,g)_H+\frac1{   \la_{\ell-1} {N_{\ell-1}}} \sum_{i=1}^{N_{\ell-1}}
\s\Big ( h_i^{(\ell-1)}(\mathbf x) \Big )
 \int_{\UU}
\s\Big ( h_i^{(\ell-1)}(\mathbf y) \Big )  g(\by) d\by,
\end{equation}
   for $\ell= 2,\dots,L+1$ and for any $g\in H$. We see that \eqref{conditional_covarianceOp_with_bias} is nothing but the translation of \eqref{conditional_covarianceOp} by the non-random linear operator $\mathbf{b}\colon g\mapsto b(1,g)_H$, defined for any given $b\geq 0$, and here used with $b=b^{(\ell)}$ at step $\ell=2,\dots,L+1$, respectively. Then the rate function  \eqref{eq:rate_function} gets modified accordingly, i.e. for any $K_1$ and $K_2$  in $\CL_1^{+,s}$,
 \[
I_{\la,b}(K_2|K_1)\colon =\sup_{ D \in \CL_\infty} \Big \{ \tr(DK_2) 
- \log\Big( \int_H e^{ \tr(D(\mathbf{b}+\frac{1}{\la}C_h))}
 \CN_H(dh|\bzero,K_1) \Big )\Big \}
\]
and so do the statements of Theorem~\ref{thm.main} and Theorem~\ref{thm.main_C0}. Given that the bias induces simply a deterministic translation of the covariances at any step of the Markov chain, this does not affect any of the proofs.

\end{remark}

\subsection{Literature review and comparison}\label{sec:literature}

In~\cite{Hanin2023} the process $\bx\mapsto h^{(L+1)}(\bx)$ is proved to converge weakly in $C^0(\UU,\RE^{D})$ to a Gaussian process with covariance function $\CK^{L+1}_{\infty}$. 
In our setting, where all weights are independent Gaussian variables (i.e., assumption~\ref{H1} holds), this corresponds to proving that the $(L+1)$-th step of the covariance process converges in $C^0(\UU^2, \RE)$ to $\CK^{L+1}_{\infty}$. Indeed, in this case, conditionally on $\CK^{L+1}_{N_{L}}$, $\bx\mapsto h^{(L+1)}(\bx)$ is a Gaussian process with covariance function $\CK^{L+1}_{N_{L}}$. Hence, our LDP in $\CC^{+,s}$ is a natural extension of this result, under our (stronger) assumptions. The assumptions in~\cite{Hanin2023} are weaker than our set of assumptions \ref{H1}--\ref{H4} in two main respects.
First, $\sigma$ is only required to satisfy polynomial bounds on its first derivatives. We expect that stronger growth conditions on $\sigma$ are necessary to obtain an LDP, although we cannot yet determine whether~\ref{H4} is optimal (in the Gaussian case).
Second, the weights in~\cite{Hanin2023} are not necessarily Gaussian (i.e.,~\ref{H1} does not hold), which means that the covariance process is no longer a Markov process. Consequently, our approach cannot be directly applied, and we expect that a more substantial extension would be required to obtain large deviation results in this case.

In  \cite{favaro2023quantitative} the same process is studied in Sobolev spaces, let us briefly discuss the differences between this choice and the one in our work. Our choice of $\CL_1^{+,s}$ equipped with the trace norm is motivated by making various tightness arguments as natural as possible. In contrast, \cite{favaro2023quantitative} has a focus on quantitative bounds and various neural network derivatives, which makes Sobolev spaces the natural choice.

To the best of our knowledge, the first work addressing large (and moderate) deviations in the context of deep neural networks is~\cite{MaPaTo24}, with an extension to the case of activation functions with linear growth in~\cite{Vog24}. These studies analyze the same neural network model considered here, focusing on large deviations of the law of the rescaled output given a finite input set $\UU_P$. Under assumption~\ref{H2} (which also applies to their setting), it is known that, conditional on the second-to-last step of the iteration (layer $L$), the network output is Gaussian. The observed deviations therefore arise from two sources: the deviations of this Gaussian distribution, and the deviations of the associated random covariance matrix from its deterministic limit (see \eqref{LLN_COV}). This structure is somewhat implicit in the form of the rate function given in~\cite[Theorem~2.1]{MaPaTo24}, which can be interpreted as the result of a contraction principle applied to the joint large deviation principle of the collection of covariance matrices and the rescaled output. Extending the analysis from the finite-dimensional setting of covariance matrices to the infinite-dimensional setting involving functional limits introduces a number of non-trivial challenges. These are addressed by carefully selecting the appropriate function space and topology. However, in the present work, we do not yet cover the case of  activation functions 
with linear growth
in this infinite-dimensional setting—a task we leave for future research.

In this regard, we note that in \cite{BassettietalJMLR} the LDP principle for the  covariance structure of a deeply linear network — that is, the case $\sigma(x) = x$ — was derived for both fully connected networks and simple convolutional architectures. In the linear case, \cite{BassettietalJMLR} provides a simpler representation of the covariance structure, which depends only on a random $D \times D$ covariance matrix, independently of the number of inputs. An explicit Markovian representation for this covariance matrix is given in the form of a product of square root random matrices with Whishart distribution. This leads to a very explicit form of the rate function, which we do not expect to be able to obtain 
in our general setting.

\section{LDP under the Posterior distribution}
\label{sec:post}

Having established the LDP for the prior covariance, we now turn to the posterior setting under a Gaussian likelihood. In Section \ref{Sec:post_A}, we derive the general form of the posterior for the random covariance functions, while in Section \ref{Sec:post_B} we establish the corresponding large deviation principles—both for the unscaled posterior and under mean-field scaling.

\subsection{The Posterior distribution of a Bayesian neural network}\label{Sec:post_A}
In a  supervised learning problem one has  a training set $\{\mathbf x_\mu ,\mathbf y_\mu\}_{\mu=1}^P$, where each $\mathbf x_\mu \in \mathbb R^{N_0 }$ 
has  the corresponding labels (response) $\mathbf y_\mu \in \RE^D$. 
The prior on $\vartheta = \{W^{(\ell-1)}_{ij}: \ell=1,\dots,L+1;  i=1,\dots, N_\ell, j=1,\dots, N_{\ell-1} \}$, say $P_{N,\mathrm{prior}}(d\theta)$,  induces a prior on all the other random elements, in particular on the 
 network outputs  $h^{(L+1)}(\mathbf x)$
  as well on  the random covariance 
 $\{\CK^{\ell}_{N_{\ell-1}}(\bx,\bx'): (\bx,\bx') \in \UU^2\}$. 
Recall that the $h^{(L+1)} (\mathbf x)$'s are functions of the parameters
$\vartheta$.

In order to perform Bayesian learning for the network parameters, one chooses a {\it likelihood for the labels given the inputs and the outputs}, in what follows denoted by $\CL(\by_{1},\dots, \by_P|\bs_{1},\dots,\bs_{P})$. In probabilistic terms, the function
\[
(\by_{1},\dots, \by_P) \mapsto \CL(\by_{1},\dots, \by_P|\bs_{1},\dots,\bs_{P})
\] represents the conditional density of the random responses 
$[\bY_1,\dots,\bY_P]$ in a generic 
point $[\by_1,\dots,\by_P]$
given the outputs  $h^{(L+1)} (\mathbf x_1)=\bs_{1},\dots,h^{(L+1)} (\mathbf x_P)=\bs_{P}$.
In analogy to a network trained with a quadratic loss function, one can consider the Gaussian likelihood
\begin{equation}\label{gaussian_likelihood}
\CL(\by_{1},\dots, \by_P|\bs_{1},\dots,\bs_P)  =\Big(\frac{\beta}{2\pi}\Big)^{{\frac{DP}2}} e^{-\frac{\beta}2\sum_{\mu=1}^P \|\bs_{\mu}-\by_{\mu}\|^2}, 
\end{equation}
with $\beta>0$.  Note that this corresponds to assuming the Gaussian error model:
\[
\bY_\mu
=h^{(L+1)}(\bx_\mu)+\bm{\varepsilon}_\mu \quad \bm{\varepsilon}_\mu \stackrel{iid}{\sim} \CN(\mathbf{0},\beta^{-1}\mathbf{1}_D) \quad \mu=1,\dots,P,
\]
 $\mathbf{1}_D$ being the $D$ dimensional identity matrix. 

The core of Bayesian learning is captured by the posterior distribution  of $\vartheta$, i.e. 
the conditional distribution of  $\vartheta$ 
given $[\bY_1,\dots,\bY_P]=[\by_{1},\dots,\by_P]$. 
The  posterior distribution 
 of  
$\vartheta$ is by Bayes theorem 
\begin{equation}\label{posterior_infty} 
P_{N,\mathrm{post}}(d\theta |\by_1,\dots,\by_P) 
:=\frac{ \CL(\by_1,\dots,\by_P| \bs_{1},\dots,\bs_{P} ) P_{N,\mathrm{prior}}(d\theta) 
}{ \int \CL(\by_1,\dots,\by_P| \bs_{1},\dots,\bs_{P}) P_{N,\mathrm{prior}}(d\theta)  } 
\end{equation}
where 
$\bs_\mu=h^{(L+1)}(\bx_\mu)\in \RE^{D}$ with $\mu=1,\dots,P$.

Here we are interested in the posterior of the empirical covariance process
\[
\boldsymbol{\mathcal{K}}_N:=(\CK^2_{N_1},\dots, \CK^{L+1}_{N_L}).
\]
In order to describe its posterior distribution  we need some more notation. 
For a continuous kernel function $\CKK$  in $\CC^{+,s}$,  define the $DP \times DP$ covariance matrix 
\begin{align}
	\label{eq:sig}
\Sigma(\CKK)=[\CKK(\bx_\mu,\bx_\nu)]_{\mu,\nu=1}^P \otimes \mathbf{1}_D
\end{align}
where $\otimes$ denotes the Kronecker product of matrices. 
Moreover, writing   $\vvec[\mathbf{A}]$  for the operation of stacking the columns of matrix $\mathbf{A}$ into a column  vector, 
 set
$\by_{1:P}=\vvec[\by_1,\dots,\by_P]$
and define the function 
\[
\Psi (\mathbf{A}|\by_{1:P})=
\beta \by_{1:P}^\top ( \eI_{DP}+\beta \mathbf{A})^{-1} \by_{1:P}
 +\log(\det(\eI_{DP} +\beta \mathbf{A}))
 \]
where $\mathbf{A}$ is a semi-positive symmetric 
$DP \times DP$ matrix.  

\begin{proposition} \label{prop:posterior}
Assuming the Gaussian likelihood described in \eqref{gaussian_likelihood}, the conditional 
 distribution of $\boldsymbol{\mathcal{K}}_N =(\CK^2_{N_1},\dots, \CK^{L+1}_{N_L})$ given $\bY_{1:P}=\by_{1:P}$ is 
\[
 \mathcal Q_N(d{\CKK}^2 \cdots d{\CKK}^{L+1} | \by_{1:P}) =% 
 \frac{e^{-\frac12 \Psi (\Sigma(\CKK^{L+1})|\by_{1:P}) }\mathcal Q_N(d{\CKK}^2 \cdots d{\CKK}^{L+1})  }{\int e^{-\frac12 \Psi (\Sigma(\CKK^{L+1})| \by_{1:P}) } \mathcal Q_N(d{\CKK}^2 \cdots d{\CKK}^{L+1}) } 
\]
where 
$\mathcal Q_N(\cdot)$ is the prior distribution for $\boldsymbol{\mathcal{K}}_N$
on $\CC^{+,s}\times \dots \times \CC^{+,s}$.
\end{proposition}

For the sake of readability, we postpone the proof of Proposition~\ref{prop:posterior} to Section~\ref{sec:proof_posterior}.

\subsection{Posterior LDP for the covariance process}\label{Sec:post_B}

Starting from the LDP for the joint law 
$\mathcal{Q}_N$, it is easy to derive the corresponding LDP for the posterior distribution of  $\boldsymbol{\CK}_N$. This requires indeed a simple adaptation of the well known large deviations result which goes under the name of Varadhan's Lemma~\cite{Var66}. See Proposition~\ref{VaradhanModified} in the appendix for the precise version of this result that we use in the following.

\begin{proposition}[Posterior LDP]
\label{thm.main_C0_posteriod}
Assume that $\UU \subset \RE^{N_0}$ is compact and  \ref{H1}-\ref{H2}-\ref{H3}-\ref{H4}, with
$r<2$ in \ref{H2}.
Then, the sequence of posterior distributions  of $\{(\CK^2_{N_1},\dots, \CK^{L+1}_{N_L})\}_N$, that is $\{\mathcal{Q}_N(\cdot| \by_{1:P})\}_N$, 
 satisfies   an LDP on 
    $\CC^{+,s} \times \dots \times \CC^{+,s}$ with speed $N$ and good rate 
function $\mathcal{I}$
given in {\rm Theorem \ref{thm.main_C0}}.
\end{proposition}

\begin{proof}
To apply Proposition~\ref{VaradhanModified}, we need to prove that $\Psi (\Sigma(\cdot)|\by_{1:P})$ is non-negative and locally bounded. First, for all $\CKK\in \CC^{+,s}$, $\Sigma(\CKK)$ is a positive semidefinite matrix, hence $\det(\eI_{DP} +\beta \Sigma(\CKK))\geq 1$ and $( \eI_{DP}+\beta \Sigma(\CKK))^{-1}$ is positive semidefinite as well, proving non-negativity of the function. Now, let $\mathcal{B}_R(\mathcal{K})$ the open ball of radius $R>0$ around  $\CKK\in \CC^{+,s}$, $\det(\eI_{DP} +\beta \Sigma(\CKK'))$ is clearly bounded on $\mathcal{B}_R(\CKK)$, while the first term is uniformly bounded since $\by_{1:P}^\top ( \eI_{DP}+\beta \Sigma(\CKK))^{-1} \by_{1:P}\leq \by_{1:P}^\top \by_{1:P}<\infty$. This implies local boundedness of $\Psi (\Sigma(\cdot)|\by_{1:P})$. Then $\Psi (\Sigma(\cdot)|\by_{1:P})$ plays the role of $\rho$ from Proposition~\ref{VaradhanModified} and the claim holds.
\end{proof}

This disappointing, though not unexpected, result shows that the LDP under the posterior remains the same as under the prior—i.e., one recovers the same rate function. This can be interpreted as yet another manifestation of the \emph{laziness} of the infinite-width asymptotic regime, as the rate function is unaffected by the training set.

Probably the easiest way to escape the lazy-training infinite-width limit is to consider the so-called \emph{mean-field parameterization}, 
see e.g. 
\cite{ChizatLazy, doi:10.1073/pnas.1806579115, GEIGER20211, Geiger_2020, yang2020feature,NEURIPS2022_d027a5c9} for networks trained using gradient descent and 
\cite{BassettietalJMLR,rubin2024a,lauditi2025adaptive} in the Bayesian setting. 
In this parameterization, 
the loss and output functions are rescaled as 
\begin{equation*}
    \CL_N(\by_{1:P}| \bs_{1:P},\beta) :=\CL(\by_{1:P}| \bs_{1:P}/\sqrt{N},N \beta)
    \label{GLikeFL}
\end{equation*}
which in the Gaussian case reads 
\begin{equation}\label{gaussian_likelihoodMF}
 \CL_N(\by_{1:P}| \bs_{1:P},\beta)   =\Big(\frac{N\beta}{2\pi}\Big)^{{\frac{DP}2}} e^{-\frac{\beta}2\sum_{\mu=1}^P \|\bs_{\mu}-\sqrt{N}\by_{\mu}\|^2}. 
\end{equation}

The mean-field parameterization exhibits some pathological behavior in the Bayesian setting.
In a sense, the scale of the prior is incorrect, since at the prior level it forces to zero both the error and the distribution of the network   in the limit. Nevertheless,  the posterior of the random covariance exhibits a well-defined and non-trivial limiting behavior. 
Comparing the large deviation asymptotics of the mean-field  posterior covariance with those in the lazy-training infinite-width limit, one recognizes the presence of additional terms, which can be interpreted as an instance of feature learning. This has been already proved for a linear network in 
\cite{BassettietalJMLR}. Here we prove a similar result for 
a general activation function $\sigma$  at the functional level.

Under the mean-field scaling, the posterior distribution  of 
$\boldsymbol{\CK}_{2:L+1}$ is 
\[
 \mathcal{Q}_N^{mf}(d{\CKK}^2 \cdots d{\CKK}^{L+1})  | \by_{1:P}) =% 
 \frac{e^{-\frac12 \Psi_N (\Sigma(\CKK^{L+1})|\by_{1:P}) }\mathcal{Q}_N(d{\CKK}^2 \cdots d{\CKK}^{L+1})   }{\int e^{-\frac12 \Psi_N (\Sigma(\CKK^{L+1})| \by_{1:P}) } \mathcal{Q}_N(d{\CKK}^2 \cdots d{\CKK}^{L+1})  } 
\]
with 
\[
\Psi_N (\mathbf{A}|\by_{1:P})=
N \beta \by_{1:P}^\top ( \eI_{DP}+\beta \mathbf{A})^{-1} \by_{1:P}
 +\log(\det(\eI_{DP} +\beta \mathbf{A})).
 \]
Due to the peculiar form of $\Psi_N$, 
the appearance of 
an additional term in the rate function is 
a consequence of the already mentioned   variant of Varadhan's Lemma (Proposition \ref{VaradhanModified} in the Appendix). 

\begin{proposition}[Posterior LDP under mean field rescaling]
\label{thm.main_C0_posteriodMF}
Assume that $\UU \subset \RE^{N_0}$ is compact and  \ref{H1}-\ref{H2}-\ref{H3}-\ref{H4}, with
$r<2$ in \ref{H2}.
Then, under the mean field 
parameterization of the likelihood \eqref{gaussian_likelihoodMF}, 
the sequence of posterior distributions $\{\mathcal{Q}_N^{mf}(\cdot| \by_{1:P})\}_N$
 satisfies   an LDP on 
    $\CC^{+,s} \times \dots \times \CC^{+,s}$ with speed $N$ and good rate 
function
\[
\mathcal{I}_{mf}(\CKK_2,\dots,\CKK_{L+1})=\mathcal{I}(\CKK_2,\dots,\CKK_{L+1})
+\beta \by_{1:P}^\top ( \eI_{DP}+\beta \Sigma(\CKK_{L+1}))^{-1} \by_{1:P}-\mathcal{I}_0
\]
for all $(\CKK_2,\dots,\CKK_{L+1})\in \CC^{+,s}\times\dots\times\CC^{+,s}$, 
where $\mathcal{I}$ is 
given in {\rm Theorem \ref{thm.main_C0}}
and 
\[
\mathcal{I}_0= 
\inf_{\CKK_2,\dots,\CKK_{L+1}} \{ I(\phi(\CKK_2),\dots,\phi(\CKK_{L+1}))
+\beta \by_{1:P}^\top ( \eI_{DP}+\beta \Sigma(\CKK_{L+1}))^{-1} \by_{1:P}\}.
\]
 \end{proposition}

\begin{proof}
We want to apply Proposition~\ref{VaradhanModified} with $\Phi_0(\CKK):=\beta \by_{1:P}^\top ( \eI_{DP}+\beta \Sigma(\CKK))^{-1} \by_{1:P}$ and $\rho(\CKK):=\log(\det(\eI_{DP} +\beta \Sigma(\CKK)))$. Following the proof of Proposition~\ref{thm.main_C0_posteriod}, we see that $\rho$ is non-negative and locally bounded. Boundedness and continuity of $\Phi_0$ follow as well. Hence, the application of Proposition~\ref{VaradhanModified} gives the claim.
\end{proof}

\section{Main ingredients for the proofs}\label{sec:main_ingredients}

The core of our approach lies in describing the vector of covariances as a Markov process with values in the most natural space, i.e. in the space of \emph{non-negative} and \emph{symmetric} trace-class  operators $\CL_1^{+,s}$. 
This formulation enables us to derive many of the limiting results in a natural way. Specifically, we apply limit theorems for sums of independent (though not identically distributed) random variables in such space, namely~\cite[Theorem~2]{Bolthausen1984}. The corresponding conditional versions of these theorems allow us to iteratively ``patch together'' the results, step by step, along the trajectory of the Markov chain.
In Section~\ref{sec:conditional} and in Section~\ref{ss:LLNLDP} we state such results under our notation and assumptions. In Section~\ref{ss:expTight} we state as well two exponential tightness results: one is  needed to iterate the LDP on $\CL_1^{+,s}$ and the second one to lift the result from $\CL_1^{+,s}\times\dots\times \CL_1^{+,s}$ to $\CC_1^{+,s}\times\dots\times \CC_1^{+,s}$. This section serves as preliminary for the proofs of the main results, that come in Section~\ref{Proofs_mainsteps}.

\subsection{Conditional limit theorems}\label{sec:conditional}

In order to get the desired LLN result for the vector of covariances, we should use a result on conditional LLN, which we state and prove here.
\begin{lemma}\label{lemmetto_conv}
Let $\{(X_{1,n},X_{2,n})\}_n$ a sequence of random vectors taking values in 
$\mathbb{X}_1 \times \mathbb{X}_2$ with $(\mathbb{X}_i,\CX_i)$ Polish. 
Let $\nu_n(x,dy)=\PP(X_{2,n} \in dy  |X_{1,n}=x)$. Assume that whenever $x_n \to x$ then 
$\nu_n(x_n,dy) \Rightarrow \delta_{G(x)}(dy)$,  
for some measurable $G:\mathbb{X}_1 \to \mathbb{X}_2$. 
As $n\to\infty$, if $X_{1,n} \to x_0$ in probability,
then $(X_{1,n},X_{2,n})$ converges in law, and hence in probability, to $(x_0,G(x_0))$. 
\end{lemma}

\begin{proof} It suffices to test the convergence for $\phi(x,y)=\phi_1(x) \phi_2(y)$ with $\phi_i$ bounded and continuous. 
Set $\Phi_{2,n}(x) :=\int_\YY \phi_2(y) \nu_n(x,dy)$. 
	Let $\{n_k\}$ be a subsequence such that $X_{n_k} \to x_0$ a.s.  Now by hypothesis if $x_n \to x$, then 
\[
\Phi_{2,n}(x_n) :=\int_\YY \phi_2(y) \nu_n(x_n,dy) \to
\phi_2(G(x)).
\]
Therefore, $\Phi_{2,n_k}(X_{n_k}) \to \phi_2(G(x_0))$ a.s.. Since 
 $\Phi_{2,n}(x) \leq \|\phi_2\|_\infty$, dominated convergence  
 gives
 \[
\E[\phi_1(X_{n_k})\phi_2(Y_{n_k})]=\E[\phi_1(X_{n_k}) \Phi_{2,n_k}(X_{n_k})] \to \E[\phi_1(x_0)\phi_2(G(x_0))].
 \]
 Since the limit is independent on the specific subsequence $n_k$ the thesis follows. 
\end{proof}

Similarly, in order to get a LDP for the joint distribution of the covariances in each layer, we shall apply a results from~\cite{Cha97} on conditional LDP.  Let us first give a useful definition.

\begin{definition}[LDP continuity condition] 
	\label{def:LDP_cont}
Let $(\BX_i, \mathcal{X}_i)_{i=1,2}$ be Polish spaces with associated Borel $\s$-algebras. A sequence of transition kernels $\{\nu_n :  \BX_1 \times  \mathcal{X}_2  \to [0, 1]\}_{n\geq 1}$ is said to satisfy the \emph{LDP continuity condition} with rate function $I_{2|1}(\cdot|\cdot)$ if:
\begin{enumerate}
    \item \label{item:good} For each $x_1 \in \BX_1$, $I_{2|1}(\cdot|x_1)$ is a good rate function on $\BX_2$.
    \item \label{item:LDP_conditional} For each $x_1 \in \BX_1$ and each sequence $x_{1,n} \to x_1$, we have that $\{\nu_n(x_{1,n},\cdot)\}_{n\geq 1}$ satisfies an LDP on $\BX_2$ with rate function $I_{2|1}(\cdot|x_1)$.
    \item \label{item:lower_semicont} The mapping $(x_1, x_2) \mapsto I_{2|1}(x_2|x_1)$ is lower semi-continuous.
\end{enumerate}
\end{definition}
Given this definition, the result from~\cite{Cha97} reads as follows.

\begin{proposition}[Theorem~2.3, \cite{Cha97}]\label{thm:joint_LDP_general}
    Let $\{\mu_n\}_n$ be a sequence of probability measures on $\BX_1$, satisfying a LDP with good rate function $I_1$. Suppose that $\{\nu_n\}_n$ satisfies the LDP continuity condition with rate function $I_{2|1}(\cdot|\cdot)$. Then:
\begin{enumerate}
    \item The sequence of measures $\{\xi_n\}_n$ defined by 
    \[
    \xi_n(A \times B) = \int_A \nu_n(x, B) \,\mu_n(dx)
    \]
    satisfies a weak LDP with rate function $I_{1,2}(x_1,x_2)= I_1(x_1) + I_{2|1}(x_2|x_1)$.
    \item If $I_{1,2}$ is a good rate function, then $\xi_n$ satisfies an LDP.
    \item The sequence of marginal measures 
    \[
    \xi_n^{(2)}(B) = \int_{\BX_1} \nu_n(x, B) \, \mu_n(dx)
    \]
    satisfies an LDP with rate function $I_2(\cdot) = \inf_{x_1} \{I_1(x_1) + I_{2|1}(\cdot|x_1)\}$.
\end{enumerate}
\end{proposition}
We refer the reader to~\cite{Cha97} for the proof.
\subsection{LLN and LDP for a triangular array  of functions of  Gaussian variables in $H$}
\label{ss:LLNLDP}

We state  a law of large numbers  for a triangular array of random elements taking values in $\CL^{+,s}_1$, tailored to our purposes. To do so, recall the operator defined in \eqref{def_Cf}, which we denote as $C_f$ for any given function $f$ on $\UU$. Under stronger assumptions, that is $r<2$ in \ref{H2}, a large deviation principle holds as well. We state it here, together with some properties of its rate function. This is an adaptation of results from~\cite{Bolthausen1984}. We postpone the proofs to Section~\ref{ss:tight}.

Given a random (measurable) element $X$ defined on 
a probability space $(\Omega,\CF,\PP)$  with values in $(\CL_1,\CBB(\CL_1))$, 
let us recall that $X$ is said to be  Bochner integrable 
if $\E[\|X\|_1]<+\infty$. Under this condition 
the expectation of $\E[X]=\int_\Omega X(\omega) \PP(d\omega)$ is a well-defined 
element of $\CL_1$, where the integral is the so-called Bochner integral, 
see, e.g., II.2 in \cite{DiestelUhl}.

\begin{proposition}\label{prop:LLN_cond}
Assume 
 \ref{H1}-\ref{H2}-\ref{H3}.  
Let 
$\{Z_{k,n}\}_{n\geq 1,k\leq n}$ be a triangular array of  random variables  taking values in 
$H$, where for all $k\leq n$, $Z_{k,n} \stackrel{ind}{\sim} \CN_H(\bzero,K_n)$ 
and $K_n \to K$ in $\CL_1^{+,s} \subset \CL_1$ as $n\to\infty$. 
Then, 
\[
S_n \colon = \frac1{n}\sum_{k=1}^n C_{Z_{k,n}}
\]
converges in probability to $\E[C_{Z_\infty}]$ when $n\to\infty$, where 
$Z_\infty  \sim \CN_H(\bzero,K)$ and $\E[C_{Z_\infty}]$ is the Bochner integral of 
$C_{Z_\infty}$. 
\end{proposition}

The next result  follows easily 
from~\cite[Theorem~2]{Bolthausen1984}, reported in Appendix as Proposition \ref{theoBolth} for completeness. 

\begin{proposition}\label{prop:LDP_cond}
Assume 
 \ref{H1}-\ref{H2}-\ref{H3}, with 
 $r<2$ in \ref{H2}. 
Let 
$\{Z_{k,n}\}_{n\geq 1,k\leq n}$ be a triangular array of  random variables  taking values in 
$H$, where for all $k\leq n$, $Z_{k,n} \stackrel{ind}{\sim} \CN_H(\bzero,K_n)$ 
and $K_n \to K$ in $\CL_1^{+,s} \subset \CL_1$ as $n\to\infty$. 
Then, for $\la>0$, the sequence of laws of 
\[
S_n \colon = \frac1{n \la}\sum_{k=1}^n C_{Z_{k,n}}
\]
 satisfies 
    a LDP   with speed $n$ and good rate 
function 
\[
I_{\lambda}(K^\prime|K)\colon =\sup_{ D \in \CL_\infty} \Big \{ \tr(DK^\prime) 
- \log\Big( \int_H e^{ \la^{-1}\tr(DC_h)}
 \CN_H(dh|\bzero,K) \Big )\Big \}
\]
for any $K^\prime \in \CL_1$.  
\end{proposition}

The rate function from the above LDP is the key ingredient for our Theorem~\ref{thm.main}, i.e. the one defined in \eqref{eq:rate_function}. It satisfies two important properties: it is equal to infinity when  $K^\prime\notin \CLL$ and it is lower-semicontinuous as a function in the two arguments $(K,K^\prime)\in \CLL\times \CLL$. We state these properties in the following lemma, see Section~\ref{ss:tight} for the proof.

\begin{lemma}\label{lem:infinity_rate_function} Fix $\lambda>0$. Let $K_1\in \CLL$, then  $I_{\lambda}(K_2|K_1)=+\infty$ for all $K_2\notin \CLL$. Moreover, 
for every sequence $\{(K^{(n)}_1,K^{(n)}_2)\}_n$ of elements in $\CLL\times \CLL$ that converges to $(K_1,K_2)\in \CLL\times \CLL$, 
 $\liminf_{n\to\infty} I_{\lambda}(K^{(n)}_2|K^{(n)}_1) \ge %h(B|\Phi(K))=
 I_{\lambda}(K_2|K_1)$. 
\end{lemma}

\subsection{Exponential tightness.}
\label{ss:expTight}
The aim of this section is to state exponential tightness for the sequence of laws we are interested in. Such results are crucial to get our desired large deviation statements. We postpone the proofs to Section~\ref{ss:expTight_L_proof} and Section~\ref{ss:expTight_C_proof}, respectively.

First, let us state a proposition on exponential tightness of the sequence of random vectors $(\bK^2_{N_1},\dots,\bK^{L+1}_{N_L})_{N \geq 1}$
in $\CL_1^{+,s}\times \cdots \times \CL_1^{+,s}$.
\begin{proposition}\label{Prop:exp.tight}  
Assume \ref{H1}-\ref{H2}-\ref{H3} with $r < 2$ in  \ref{H2},  
then the sequence of laws of
\[
\{(\bK^2_{N_1}, \bK^3_{N_2},\dots,\bK^{L+1}_{N_L})\}_{N}
\]
	is exponentially tight, i.e., for any $R>0$, there exists a compact set $\mathcal{B}_R\subset (\CLL)^L$ such that
\[
\limsup_{N\to\infty}\frac 1N\log \PP\left((\bK^2_{N_1}, \bK^3_{N_2},\dots,\bK^{L+1}_{N_L})\notin \mathcal{B}_R \right)\leq -R.
\] 
\end{proposition}

Secondly, we show that also the sequence of laws of  $(\CK^2_{N_1},\dots,\CK^{L+1}_{N_L})_{N \geq 1}$
is  exponentially tight in $\CC^{+,s}\times \cdots \times \CC^{+,s}$.
For this result, besides previous assumptions, we assume additionally   \ref{H4}.

\begin{proposition}\label{prop:exp_tight_C0} Assume that $\UU \subset \RE^{N_0}$ is compact and  \ref{H1}-\ref{H2}-\ref{H3}-\ref{H4}, with
$r \leq 2$ in \ref{H2}.
Then, the sequence of laws of $\{(\CK^2_{N_1},\dots, \CK^{L+1}_{N_L})\}_N$
is exponentially tight, that is for any $M<\infty$, there exists a compact set $\CC_M\subset (\CC^{+,s})^L$ such that
\[\limsup_{N\to\infty}\frac 1N\log \PP\left((\CK^2_{N_1},\dots, \CK^{L+1}_{N_L})\notin \CC_M \right)\le -M.
\]
\end{proposition}

\section{Proofs of main limiting results}
\label{Proofs_mainsteps}

In this section we gather the proofs of our main results: Proposition \ref{prop:LLN_traceclass}, Theorem \ref{thm.main} and Theorem \ref{thm.main_C0}. The proofs use the fact that the transition rules~\eqref{eq:K_in_C0}  and \eqref{transitionK} define Markov chains with values in  $\CC^{+,s}$ and in $\CL_1^{+,s}$, respectively. This, together with the results from Section~\ref{sec:main_ingredients}, is crucial to get the desired results. In the process we also use some additional results, which we collect in Section~\ref{sec:proofs}.

\subsection{Proof of Proposition \ref{prop:LLN_traceclass}}\label{proof_temp_LLN_L1}

Given the Markov chain structure of $(\bK^2_{N_1},\dots,\bK^{L+1}_{N_L})$, see~\eqref{transitionK}, we will apply Lemma 
\ref{lemmetto_conv} iteratively in order to prove  Proposition \ref{prop:LLN_traceclass}. 
The Polish spaces are the spaces of \emph{self-adjoint, non-negative} trace-class operators, i.e.\ $\BX_i=(\CLL,\|\cdot\|_1)$ for $i=1,2,\dots,L$. Note that from~\ref{Fact0}, $(\CL_1,\|\cdot\|_1)$ is indeed Polish
and hence  $\CLL \subset \CL_1$ is Polish being 
a closed subset of a Polish space.

 By a direct application of Proposition~\ref{prop:LLN_cond}, the sequence  $(\bK^2_{N_1})_{N_1}$ converges in 
 probability to $\bK^2_{\infty}=\E[C_{Z^2_\infty}]$
where $Z^2_\infty \sim \CN_H(\bzero,
\bK^{1}_{\infty} )$. 
Indeed, from~\eqref{conditional_h} 
 \begin{equation}\label{K2-sum}
 \bK^2_{N_{1}} =\frac1{   \la_{1} {N_{1}}} \sum_{i=1}^{N_{1}} C_{ h_i^{(1)}},
 \end{equation}
with $\{h^{(1)}_i\}_{i=1,\dots, N_1}$  i.i.d. random variables with law  $\CN_H(\bzero,\bK^1_{N_0})$. Now, for any $K\in \CLL$,  let $\SSS_{\la_2, N_2}(\cdot|K)$ be a transition kernel corresponding to the  conditional law of 
$\bK^3_{N_2}$ given $\bK^2_{N_1}=K$, see next Lemma~\ref{regualr_CondProb_kernel}.
Using again Proposition~\ref{prop:LLN_cond}
also the second assumption in Lemma \ref{lemmetto_conv} is satisfied
for $G(K)=\E[C_{Z_K}]$ with $Z_K\sim \CN_H(\bzero,K)$, $K \in \CL^{+,s}$. Applying this iteratively 
one obtains 
\[
(\bK^2_{N_1},\dots,\bK^{L+1}_{N_L}) \stackrel{\PP}{\to}  (\bK^2_{\infty},\dots,\bK^{L+1}_{\infty}) \quad \text{as $N \to \infty$}
\]
where $(\bK^2_{\infty},\dots,\bK^{L+1}_{\infty})$ are recursively defined as
\[
\bK^\ell_{\infty}=\E[C_{Z^{\ell}_\infty}]
\qquad Z^{\ell}_\infty \sim \CN_H(\bzero,
\bK^{\ell-1}_{\infty} ) \quad \ell=2,\dots,L+1.
\]
It remains to check that 
\[
\bK^\ell_{\infty}=\phi(\CK^\ell_{\infty}).
\]
Now $\bK^1_{\infty}=\bK^\ell_{N_0}=\phi(\CK_{N_1}^1)$ by definition, where $\phi$ is defined in \eqref{mapKerneltoCov}.
Assume that $\bK^{\ell-1}_{\infty}=\phi(\CK^{\ell-1}_{\infty})$
with $\CK^{\ell-1}_{\infty} \in \CC^{+,s}$
and note that for $g_1,g_2 \in H$
\[
\langle \bK^\ell_{\infty}g_1,g_2 \rangle
=\langle \E[C_{Z^{\ell}_\infty}] g_1,g_2 \rangle
=\E[ \langle C_{Z^{\ell}_\infty} g_1,g_2 \rangle] 
\]
where the second equality follows 
by linearity of Bochner integral since
$\bK \mapsto \langle \bK g_1,g_2 \rangle$
is a bounded linear  operator. 
Now, 
\[
\langle C_{Z^{\ell}_\infty} g_1,g_2 \rangle
=\int \sigma(Z^{\ell}_\infty(\bx_1) ) \int \sigma(Z^{\ell}_\infty(\bx_2)) g_1(\bx_2)d\bx_2 g_2(\bx_1)d\bx_1.
\]
for $Z^{\ell}_\infty \sim \CN_H(\bzero,
\bK^{\ell-1})$. Since 
 $\bK^{\ell-1}_{\infty}=\phi(\CK^{\ell-1}_{\infty})$,  
 $Z^{\ell}_\infty$ can be identified 
 (in law) with  $Z^{\ell}_\infty \sim GP(\bzero,\CK^{\ell-1}_{\infty})$
and, by Fubini theorem, 
\[
\begin{split}
\E[ \langle C_{Z^{\ell}_\infty} g_1,g_2 \rangle] 
& = \int \int \E[\sigma(Z^{\ell}_\infty(\bx_1) )\sigma(Z^{\ell}_\infty(\bx_2)) ]g_1(\bx_2) g_2(\bx_1)d\bx_2d\bx_1\\
& = 
 \int \int \CK^{\ell-1}_{\infty} (\bx_1,\bx_2) g_1(\bx_2) g_2(\bx_1)d\bx_2d\bx_1. \\
\end{split}
\]
This shows that $\bK^{\ell}_{\infty}=\phi(\CK^{\ell}_{\infty})$ and ends the proof of Proposition~\ref{prop:LLN_traceclass}.

\subsection{LDP in $\CLL\times \dots\times \CLL$: proof of Theorem~\ref{thm.main}}
\label{ss:cldp}

Given the Markov chain structure of $(\bK^2_{N_1},\dots,\bK^{L+1}_{N_L})$, see~\eqref{transitionK}, we will apply results from Section~\ref{sec:main_ingredients} iteratively in order to prove an LDP for the sequence of its law on $(\CLL)^L$. First, we prove an LDP for the law of $(\bK^2_{N_1},\bK^3_{N_2})$ on $(\CLL)^2$. Then, in general, we prove that if we have an LDP for the law of $(\bK^2_{N_1},\dots, \bK^\ell_{N_{\ell-1}})$ on $(\CLL)^{\ell-1}$, then we have an LDP for the law of $(\bK^2_{N_1},\dots, \bK^{\ell+1}_{N_\ell})$ on $(\CLL)^{\ell}$.

As before, we take  $\BX_i=(\CLL,\|\cdot\|_1)$ for $i=1,2,\dots,L$. 
First, let us prove that the sequences of laws of $(\bK^2_{N_1},\bK^3_{N_2})$ satisfies a LDP on $(\CLL)^2$. In order to apply Proposition~\ref{thm:joint_LDP_general} to our case, we should prove that the assumptions hold. Starting from 
\eqref{K2-sum}, by a direct application of Proposition~\ref{prop:LDP_cond}, the sequence of laws of $(\bK^2_{N_1})_{N_1}$ satisfies an LDP in $\CLL$ with speed $N_1$ and good rate function $I_{\la_1}(\cdot|\bK^1_{N_0})$, as defined in~\eqref{eq:rate_function}. 
 Now, for any $K \in \CLL$,  $\bK^3_{N_2}$ given $\bK^2_{N_1}=K$ has law $\SSS_{\la_2, N_2}(\cdot|K)$ (this is defined later in Lemma~\ref{regualr_CondProb_kernel}). We claim that the sequence of kernels $\{\SSS_{\la_2, N_2}(\cdot|\cdot)\}_{N_2}$ satisfies the \emph{LDP continuity condition} from Definition~\ref{def:LDP_cont} with rate function $I_{\la_2}(\cdot|\cdot)$. Let us check that all points of Definition~\ref{def:LDP_cont} are indeed satisfied:
\begin{itemize}
    \item points \ref{item:good} and \ref{item:LDP_conditional} are consequences of Proposition~\ref{prop:LDP_cond}, 
    \item point \ref{item:lower_semicont} is a consequence of Lemma~\ref{lem:infinity_rate_function}. 
    \end{itemize}
    Then, as a consequence of Proposition~\ref{thm:joint_LDP_general}, the sequence of laws of $\{(\bK^2_{N_1},\bK^3_{N_2})\}_N$ satisfies a weak LDP with speed $N$ and rate function 
\[
I(K_2,K_3)=m_1 I_{\la_1}(K_2|\bK^1_{N_0})+m_2I_{\la_2}(K_3|K_2), \quad \forall (K_2,K_3)\in \CLL\times  \CLL,
\]
since $\lim_{N\to\infty}{N_\ell}/{N}=m_\ell$, for $\ell=1,2$ because of~\ref{H3}.  
Thanks to Proposition~\ref{Prop:exp.tight},  $\{(\bK^2_{N_1},\bK^3_{N_2})\}_N$
is exponentially tight 
and hence, by
\cite[Lemma 1.2.18]{DemboZeitouni}, 
the rate function $I(K_2,K_3)$ is good and $\{(\bK^2_{N_1},\bK^3_{N_2})\}_N$ satisfies an LDP with speed $N$ and rate function $I(K_2,K_3)$.

Now, for any $\ell\leq L$, suppose that the sequence of laws of  $\{(\bK^2_{N_1},\dots, \bK^\ell_{N_{\ell-1}})\}_N$ satisfies an LDP  on $(\CLL)^{\ell-1}$ with speed $N$ and good rate function $I(K_2,\dots,K_{\ell-1})$. We shall prove that the sequence of laws of $\bK^{\ell+1}_{N_\ell}$ given $(\bK^2_{N_1},\dots, \bK^\ell_{N_{\ell-1}})$ satisfy the \emph{LDP continuity condition} from Definition~\ref{def:LDP_cont} with speed $N_{\ell-1}$ and rate function $I_{\la_{\ell-1}}(\cdot|\cdot)$. By Markov property, the law of $\bK^{\ell+1}_{N_\ell}$ given $(\bK^2_{N_1},\dots, \bK^\ell_{N_{\ell-1}})=(K_2,\dots,K_\ell)$ equals $\SSS_{\la_\ell, N_\ell}(\cdot|K_\ell)$ for any $(K_2,\dots,K_\ell)\in (\CLL)^{\ell-1}$. Following the steps above, we see that the sequence of kernels $\{\SSS_{\la_\ell, N_\ell}(\cdot|\cdot)\}_{N_\ell}$ satisfies the \emph{LDP continuity condition} from Definition~\ref{def:LDP_cont} with speed $N_\ell$ and rate function $I_{\la_\ell}(K_{\ell+1}|K_\ell)$ for any $\left((K_2,\dots,K_\ell),K_{\ell+1}\right)$ in $(\CLL)^{\ell-1}\times \CLL$. Then, recalling~\ref{H3}, we apply again Propositions~\ref{thm:joint_LDP_general} and \ref{Prop:exp.tight} to prove that the sequence of laws of $\{(\bK^2_{N_1},\dots, \bK^{\ell+1}_{N_\ell})\}_N$ satisfies an LDP  on $(\CLL)^{\ell}$ with speed $N$ and the desired rate function, which is good. This finishes the proof of Theorem~\ref{thm.main}.

\subsection{LDP in $\CC^{+,s}\times \dots\times \CC^{+,s}$: proof of Theorem~\ref{thm.main_C0}}

The proof of Theorem~\ref{thm.main_C0} follows directly by the \emph{inverse contraction principle}~\cite[Theorem~4.2.4]{DemboZeitouni}. Indeed, let us call $\phi_L\colon (\CC^{+,s})^L\to (\CLL)^L$ the function defined as \[\phi_L({\CKK}_2,\dots,{\CKK}_{L+1})\colon =\left(\phi({\CKK}_2),\dots,\phi({\CKK}_{L+1})\right),\qquad \forall \, ({\CKK}_2,\dots,{\CKK}_{L+1})\in \CC^{+,s}\times \dots\times \CC^{+,s}.\]  This function is a continuous injection as a consequence of Lemma~\ref{continuity_Phi}. Moreover, by Proposition~\ref{prop:exp_tight_C0}, the sequence of laws of $(\CK^2_{N_1},\dots, \CK^{L+1}_{N_L})$ is exponentially tight and  $\phi_L(\CK^2_{N_1},\dots, \CK^{L+1}_{N_L})$ has the same distribution of $(\bK^2_{N_1},\dots, \bK^{L+1}_{N_L})$. Since the sequence of laws of $(\bK^2_{N_1},\dots,\bK^{L+1}_N)$
	satisfies an LDP on $(\CLL)^L$ with speed $N$ and rate function $I(\cdot)$ by Theorem~\ref{thm.main}, the \emph{inverse contraction principle} implies that the sequence of laws of $(\CK^2_{N_1},\dots, \CK^{L+1}_{N_L})$ satisfies an LDP on $\CC^{+,s}\times \dots\times \CC^{+,s}$ with speed $N$ and rate function $\mathcal{I}(\cdot)=I(\phi_L(\cdot))$. This ends the proof of Theorem~\ref{thm.main_C0}.

\section{Additional results and proofs}
\label{sec:proofs}
 After having discussed the main ideas in Section \ref{Proofs_mainsteps}, we now give the remaining details. First, in Section \ref{sec:pprel}, we discuss a couple of preliminary results and general tightness conditions for Gaussian process. These conditions are useful for discussing LLN and LDP conditions for Gaussian processes from Section~\ref{ss:LLNLDP}. In Section \ref{ss:tight} we give the proofs of the results stated in Section~\ref{ss:LLNLDP}.  In Sections \ref{ss:expTight_L_proof} and \ref{ss:expTight_C_proof}, we prove the exponential tightness results in $\CLL\times \dots \times \CLL$  and $\CC^{+,s}\times \dots\times \CC^{+,s}$, respectively. These results  are stated in Section~\ref{ss:expTight}. Finally, in Section~\ref{sec:proof_posterior}, we give the proof of Proposition~\ref{prop:posterior}.

\subsection{ Preliminary results}
\label{sec:pprel}

\begin{lemma}\label{Lemma_cont_Cf}
Assume \ref{H2} with $r\le 2$.  
The function  $f \mapsto C_f$, defined in~\eqref{def_Cf}, defines a 
continuous application from $(H,\|\cdot\|_{H})$ to $(\CL_1^{+,s},\|\cdot\|_1)$.  Moreover, 
\begin{equation}\label{ineq_s(h)}
\|C_f\|_1=\|\s(f)\|_H^2\leq   A\Big (1+ \|f\|_H^{r} \Big ).
\end{equation}
\end{lemma}

\begin{proof}
Note that if \ref{H2} holds for some $r \le 2$,  it holds as well  
for $r=2$ (with a different $A$), so let us assume that $r=2$. 
If $f \in H$,  then  by \ref{H2} (with $r=2$) one has
$\int_\UU  \s (f(y))^2 dy \le A \int_\UU (1+|f(y)|^2) dy <+\infty$ 
since   $f$ is in $H$. This shows that  $\sigma(f):=\s \circ f$  is in $H$ and 
$(\sigma(f) ,g)_H $ is well-defined. 
Moreover,  $C_f$ is clearly positive since
$(C_fg,g)_H=(\sigma(f) ,g)_H^2\ge 0$. 
The symmetry follows by noticing that 
$(C_fg_1,g_2)_H=(C_fg_2,g_1)_H$ for every $g_1,g_2 \in H$.  
We now see that it is also trace-class.  Let $(\ee_i)_i$ be an orthonormal basis for $H$ and write
\begin{equation}\label{traceCf}
\tr(C_f)=\sum_i (C_f \ee_i,\ee_i)_H
=\sum_i (\sigma(f) ,\ee_i)_H^2
=\|\s(f)\|_H^2 <+\infty
\end{equation}
where the last equality is Parseval's identity.  In conclusion, 
     $f \mapsto C_f$ maps $H$ into $\CL_1^{+,s}$.  
     Assume that $f_n \to f$ in $H$. 
     Since we are assuming that $\UU$ is bounded, 
$f_n \to f$ in $L^2(\UU)$ yields that $f_n \to f$ in measure, so that 
$\s(f_n) \to \sigma(f)$  in measure since $\s$ is continuous. 
Using again  $\s(f_n(x))^2\le A(1+|f_n(x)|^2)$ and $f_n \to f$ in $L^2(\UU)$, 
	generalized dominated convergence theorem yields that $\int_\UU \s(f_n(x))^2dx \to \int_\UU \s(f(x))^2dx $,
and hence also $\s(f_n) \to \s(f)$  in $H$. 
Setting $\CC_f(x,y)=\s(f(x))\s(f(y))$, one has $\CC_f \in L^2(\UU^2)$ and the
corresponding Hilbert-Schmidt on $H$ is the 
operator   $C_f=\phi(\CC_f)$. Similarly $C_{f_n}-C_f$ is canonically identified 
with the kernel $\CC_{f_n}-\CC_f$. Using the isometry 
of the Hilbert-Schmidt operators on $H=L^2(\UU)$ with  $L^2(\UU^2)$, see \ref{item:Hilbert-Schmidt} in the Appendix, 
one has 
     $\|C_{f_n}-C_f\|_2=
     \|\CC_{f_n}-\CC_f\|_{L^2(\UU^2)}$, where $\| \cdot\|_2$ is the Hilbert-Schmidt norm. 
    At this stage using the fact that $\s(f_n) \to \s(f)$ in $L^2(\UU)$, it is easily checked that $  \|\CC_{f_n}-\CC_f\|_{L^2(\UU^2)} \to 0$, since $\|\CC_{f_n}-\CC_f\|_{L^2(\UU^2)} \leq C  \|\sigma(f_n)-\sigma(f)\|_{L^2(\UU)} \sup_n\|\sigma(f_n)\|_{L^2(\UU)}$ for some constant $C$.
    By \eqref{traceCf} also $|\tr(C_{f_n})-\tr(C_f)|
     = \big | \|\s (f_n)\|_H^2-\|\s(f)\|_H^2 \big | \to 0$, and  Lemma \ref{Lemma_eqvivalent_conv} gives 
     $\| C_{f_n} -C_f \|_1 \to 0$. 
 Using positivity and \eqref{traceCf}, one gets
$\|C_f\|_1=\tr(C_f)=\|\s(f)\|_H^2$.
Then,  by \ref{H2}
\begin{equation*}
\|\s(f)\|_H^2
=\int_\UU \s(h(x))^2dx 
\le A\Big (1+\int_\UU |f(x))|^{r} dx \Big )
=  A\Big (1+ \|f\|_H^{r} \Big ).
\end{equation*}
\end{proof}

For a Polish space $\BX$, let $\mathcal{M}_1(\BX)$ be the space of all the probability measures on the Borel $\s$-field $\CBB(\BX)$, endowed with the topology of the weak convergence.

Let us recall or introduce the following maps:
\begin{enumerate}[label=(M\arabic*)]
\item \label{def:phi}$\phi \colon \CC^{+,s} \subset C^0(\UU^2,\RE) \to \CL_1^{+,s}$ defined by
$\phi(\CKK)\colon =K$, see 
\eqref{mapKerneltoCov};
    \item \label{def:gamma} 
$\gamma: \CLL \to \mathcal{M}_1(H)$ defined by  $ \gamma(K)\colon = \CN_H(\mathbf{0},K)$ for all $ K\in \CLL$;
\item \label{def:PHI}  $\Phi: \CLL\to \mathcal{M}_1(\CLL)$  
 defined by  $\Phi(K)=\gamma(K) \#  {C}_\cdot=\text{law}({C}_Z)$, where $Z \sim \CN_H(\mathbf{0},K)$ for all $ K\in \CLL$. 
\end{enumerate}

\begin{lemma}\label{continuity_Phi}
The  maps  $\phi$, $\gamma$ and $\Phi$ defined in \ref{def:phi}, \ref{def:gamma} and \ref{def:PHI}, respectively,
 are continuous. 
\end{lemma}

\begin{proof}
 The proof of the continuity of $\phi$ is very similar to the proof of Lemma \ref{Lemma_cont_Cf}. If $\CKK_n \in \CC^{+,s}$ converges to $\CKK$ in $C^0(\UU^2,\RE)$, then $\CKK \in  \CC^{+,s}$, moreover, since $\UU$ is compact
it follows easily that $\phi(\CKK_n) \to \phi(\CKK)$ in $\CL_2$.  
By Mercer's theorem (see  \ref{Mercer} in Appendix) $\tr(\phi(\CKK_n))=
\int_{\UU^2} \CKK_n(x,y) dx dy$, so that $\tr(\phi(\CKK_n)) \to \tr(\phi(\CKK))$. Then,
	Lemma \ref{Lemma_eqvivalent_conv} gives $\|\phi(\CKK_n)-\phi(\CKK)\|_1 \to 0$.
By  \ref{Fact1} if $K_n \to K$ in $\CLL$ then
 $\CN_H(\mathbf{0},K_n)$ converges weakly 
to $\CN_H(\mathbf{0},K)$, showing that $\gamma$ is continuous as well. 
Finally, recall that  
$\Phi(K)=\gamma(K) \#  {C}_\cdot$. Having proved that  $K\mapsto \gamma(K)=\CN_H(\mathbf{0},K)$
is continuous, since also  $f\mapsto {C}_f$ 
is  continuous (see Lemma  \ref{Lemma_cont_Cf}), the continuity of   $\Phi$ follows. 
\end{proof}

Given $\la>0$ and $N \ge 1$,
 let $\HH_{\la,N}:  \mathcal{M}_1( \CL_1^{+,s})  \to \mathcal{M}_1( \CL_1^{+,s})$
be defined by  
 \[
 \HH_{\la,N}(\nu)= \text{law} \Big (\frac1{\la N} \sum_{i=1}^N \bK_i \Big)  \qquad  \bK_i \stackrel{iid}{\sim} \nu
 \]
 and  set
 \[\SSS_{\la,N} =  :\HH_{\la,N} \circ \Phi.\]
We put $\HH_N \colon = \mathfrak H_{1, N}$ and $\SSS_N \colon = \SSS_{1, N}$ and without loss of generality we prove the following lemma for $\la=1$.

\begin{lemma}\label{regualr_CondProb_kernel}  Under the previous assumptions and notations, write  
\[
\SSS_{ N}(A|K)\colon =\SSS_N(K)(A)   \qquad  \forall A \in \CBB(\CL_1^{+,s}) \quad  \forall  K \in \CL_1^{+,s}. 
\] 
Then, $\SSS_N(\cdot|\cdot)$ 
is a probability kernel on $\CBB(\CL_1^{+,s}) \times\CL_1^{+,s}$.  
Moreover, 
for every $A \in \CBB(\CL_1^{+,s})$ and $K \in \CL_1^{+,s}$ one has 
\[
\SSS_N(A|K)=\PP \Big ( \frac1{N} \sum_{i=1}^N C_{Z_i} \in A  \Big ) \qquad Z_i \stackrel{iid}{\sim} \CN_H(\mathbf{0},K). 
\]
\end{lemma}

\begin{proof}  Since $\HH_N$ is measurable and $\Phi$ is continuous also 
the map $\SSS_N=\HH_N \circ \Phi:   \CL_1^{+,s} \to \mathcal{M}_1( \CL_1^{+,s})$  is measurable.
Under the  assumptions  $\CBB(\mathcal{M}_1( \CL_1^{+,s}))$ coincides with the smallest $\s$-field which contains all the 
evaluation map $B \mapsto p(B)$ for $B \in \CBB( \CL_1^{+,s})$. 
Since the evaluation map is measurable form $ \mathcal{M}_1( \CL_1^{+,s})$ into $[0,1]$, 
one also obtains that  $K \mapsto \SSS_N(A|K)$ is measurable, proving that $\SSS_N(\cdot|\cdot)$ is a kernel. 
The second part follows directly from  \eqref{transitionK}. 
\end{proof}
The result above holds equivalently with $\la\in (0,\infty)$ for $\SSS_{\la, N}$.
In particular, when $\la=\la_\ell$, $\SSS_{\la_\ell, N_{\ell}}(\cdot|K)$ is  the  conditional distribution of 
$\bK_{N_{\ell}}^{\ell+1}$ given $\bK_{N_{\ell-1}}^\ell=K$.

Finally, we collect here two useful facts about Gaussian measures on $H$.

\begin{lemma}\label{Lemma:compactness}
Let $\mathcal{A}$ be a compact set in $\CL_1^{+,s}$, then
$\mathscr{A}: = \{ \CN_H(\bzero,K) : K \in \CA\} \subset \mathcal{M}_1(H)$
is a tight family. Moreover, there is $t_0>0$ such that 
\[
\sup_{\gamma \in \mathscr{A}} \int_H e^{t_0 \|h\|_H^2 } \gamma(dh)=\sup_{K \in \CA} \int_H e^{t_0 \|h\|_H^2 }  \CN_H(dh |\bzero,K) <+\infty.
\]
Finally, for  every $\epsilon \in (0,2)$ and every $t>0$,
\begin{equation}\label{eq:assumption_nu_n2}
\sup_{\gamma \in \mathscr{A}}  \int_H e^{t \|h\|_H^{2-\epsilon} } \gamma(dh) <+\infty.
\end{equation}
\end{lemma}

\begin{proof}
Recall that  if $K_n \to K$ in $\CA$ (wrt the $\CL_1$ norm), 
then $\CN_H(\bzero,K_n)$ converges weakly to $\CN_H(\bzero,K_n)$, by  \ref{Fact1}. 
This means that for any sequence  $\gamma_n=\CN_H(\bzero,K_n)$ in $\mathscr{A}$, there is a subsequence 
$K_{n_k} \to K$ in $\CL_1$ and hence $\gamma_{n_k}$ is 
weakly convergent.  So that  the first claim follows from  Prohorov's theorem. 
The second part is a consequence of the Fernique’s theorem.  To see this one can apply  Theorem 3.8.11 in \cite{Bogachev}. Even if this  
theorem is stated for a sequence and not for a family,  the proof can be easily adapted. Alternatively,  
one  can apply  Theorem 2.2 in  \cite{Baldi2020},
 whose assumptions 
 are satisfied since 
 compact sets are bounded in  $H$ and hence, being $\mathscr{A}$ tight, there is $s$ such that $\gamma(\|h\|_H >s ) \le \beta<1$ with $\beta<1/2$ for all $\gamma \in \mathscr{A}$. Then, 
	\[
    \begin{split}
    \sup_{\gamma \in \mathscr{A}} \int_H e^{t \|h\|_H^{2-\epsilon} } \gamma(dh) 
    & = \sup_{\gamma \in \mathscr{A}}  \Big \{\int_H  e^{t \|h \|_H^{2-\epsilon} } \II \Big( \frac t{t_0} < \|h\|_H^\epsilon\Big) \gamma(dh) 
    +  \int_H e^{t \|h\|_H^{2-\epsilon} }
\II \Big (\frac t{t_0} \ge \|h\|_H^\epsilon \Big) \gamma(dh) \Big \}  \\
	    &    \le \sup_{\gamma \in \mathscr{A}}  \int_H e^{t_0 \|h \|_H^2 } \gamma(dh)
    +e^{t (t/t_0)^{(2-\epsilon)/\epsilon}}<+\infty,\\
    \end{split}
    \]
thereby proving \eqref{eq:assumption_nu_n2}.
\end{proof}

The previous result translates immediately in an uniform 
bound on exponential moments of $C_h$. 
Recall that, given $h \in H$,
$C_h$ is characterized by $[C_h g](x)=(\s( h ) ,g)_H \s(h (x))$ for every $g \in H$. 

\begin{lemma}\label{Lemma:compactness2}
Let $\mathscr{A}  = \{ \CN_H(\bzero,K) : K \in \CA\}$ where $\CA$  is a compact set in $\CL_1^{+,s}$
If  \ref{H2} holds with $r<2$,  
then for every $t>0$
\[
\sup_{\gamma \in \mathscr{A}}   \int_H  e^{t \| C_h  \|_1} \gamma(dh) <+\infty.
\] 
\end{lemma}

\begin{proof}
Since $r<2$ the thesis follows 
combining \eqref{ineq_s(h)} in Lemma 
\ref{Lemma_cont_Cf}  and \eqref{eq:assumption_nu_n2} in Lemma \ref{Lemma:compactness}.
\end{proof}

%
%TIGHT
%
\subsection{LLN and LDP: proofs from Section~\ref{ss:LLNLDP}}
\label{ss:tight}

We are ready now to prove the results stated in Section~\ref{ss:LLNLDP}.
We start with the proof of the LLN, which holds under assumptions  \ref{H1}-\ref{H2}-\ref{H3}.  

\begin{proof}[Proof of Proposition~\ref{prop:LLN_cond}]
Since $K_n \to K$, by \ref{Fact1} one has that   $Z_{1,n}$ converges in law 
in $H$ to $Z_{\infty} \sim \CN(0,K)$.
By Lemma~\ref{Lemma_cont_Cf} and the  continuous mapping theorem,  also $C_{Z_{1,n}}$ converges  in law 
in $\CL_1$ to $C_{Z_{\infty}}$.  
 Hence, the thesis follows from  Proposition \ref{LLNarray} 
if we prove that 
\begin{equation}\label{Bolt4}
\sup_n \E[\|C_{Z_{1,n}}\|_1^{p}]<+\infty
\end{equation} 
for some  $p > 0$.
Now, 
since $Z_{1,n} \stackrel{\CL}{\to} Z_{\infty}$ then 
$\sup_n\E[\|Z_{1,n}\|_H^{rp}]<+\infty$,   
see Thm. 3.8.11 in \cite{Bogachev}. 
Hence, using 
 \eqref{ineq_s(h)} of Lemma 
\ref{Lemma_cont_Cf},
\[
\sup_n \E[\|C_{Z_{1,n}}\|_1^{p}]
\leq A^p \sup_n \E[(1+\|Z_{1,n}\|_H^r)^p]  <+\infty.
\]
 \end{proof}

 Under stronger hypothesis, namely $r<2$ in \ref{H2},
 a LDP holds as well, which we prove here. As mentioned, this is an adaptation of~\cite[Theorem~2]{Bolthausen1984}.

\begin{proof}[Proof of Proposition~\ref{prop:LDP_cond}]
In the proof of 
Proposition \ref{prop:LLN_cond} we have shown that  $\la^{-1}C_{Z_{1,n}} \stackrel{\CL}{\to}\la^{-1} C_{Z_\infty}$.
Hence, the thesis follows from  Proposition \ref{theoBolth} in the Appendix provided that  
\begin{equation}\label{Bolt2}
\sup_n \E[e^{t\la^{-1} \|C_{Z_{1,n}}\|_1}]<+\infty
\end{equation} 
for every $t > 0$.
Lemma \ref{Lemma:compactness2}
applied to the tight family $\mathcal{A}=\{K_n : n \ge 1\}$, gives~\eqref{Bolt2}. 
\end{proof}

Finally, we prove some important properties for the rate function defined in~\eqref{eq:rate_function}.

\begin{proof}[Proof of Lemma~\ref{lem:infinity_rate_function}]
Without loss of generality we fix $\lambda=1$ and we denote with $I(K_2|K_1)$ the function $I_{\lambda}(K_2|K_1)$ in this case.
Since  
   $\P( S_n \in \CL_1^{+,s})=1$ and $\CL_1^{+,s}$ 
   is closed, the first part is
   a consequence of Lemma 4.1.5 (b) 
   \cite{DemboZeitouni}.  
   As for the second part of the statement, 
 note that  the rate function $I(K_2|K_1)$ can be written as $h(K_2|\Phi(K_1))$, where $h$ is the rate function appearing in Proposition~\ref{theoBolth} in the Appendix, 
   and $\Phi$ is the continuous function given in Lemma \ref{continuity_Phi}. By continuity,
  $\Phi(K^{(n)}_1)\to\Phi(K_1)$ in $\mathcal{M}_1(\CLL)$. Moreover 
 Lemma \ref{Lemma:compactness2} applied to $\CA=\{K^{(n)}_1 : n \ge 1\}$ gives that 
    the sequence of measures $\Phi(K^{(n)}_1)$ satisfies \eqref{Bolt2}. Hence, we can apply \cite[Lemma 1]{Bolthausen1984} which gives 
    \[
    \liminf_{n\to\infty} I(K^{(n)}_2|K^{(n)}_1)=\liminf_{n\to\infty} h\big(K^{(n)}_2|\Phi(K^{(n)}_1)\big)\ge h\big(K_2|\Phi(K_1)\big)=I(K_2|K_1).
    \]
\end{proof}

\subsection{Exponential tightness of the laws in $\CL_1^{+,s}\times \cdots \times \CL_1^{+,s}$: proof of Proposition~\ref{Prop:exp.tight}}
\label{ss:expTight_L_proof}

The aim of this section is to prove 
that the sequence of laws of the random vectors $\{(\bK^2_{N_1},\dots,\bK^{L+1}_{N_L})\}_{N \geq 1}$
is  exponentially tight in $\CL_1^{+,s}\times \cdots \times \CL_1^{+,s}$.

To prove the exponential tightness
we shall take
advantage of the following result which is a particular case of \cite{deAcosta85}. %\CH{Put in appendix or not?}

\begin{proposition}[Theorem 3.1 \cite{deAcosta85}]\label{thm:deAcosta} Assume that $\mathscr{A}$ is a tight family of probability measure over a Banach space $(E,\|\cdot\|_E)$. 
If for every $t>0$  
\[
\sup_{\mu \in \mathscr{A}} \int_E e^{t \| e  \|} \mu(de) <+\infty, 
\]
then there is a compact, convex, well balanced  set $V$ in $E$ such 
\[
\sup_{\mu \in \CA} \int_E e^{ q_V (e) } \mu(de) <+\infty
\]
for $q_V(e)=\inf\{ t \ge 0: e \in t V\}$. Under these assumptions $q_V$ is subadditive and positively homogeneous.
\end{proposition}
We use the above proposition to prove the following crucial lemma along the lines of \cite{Baldi2020}.
\begin{lemma}\label{Lemma:compactness3}
Let $\mathscr{A}  = \{ \CN_H(\bzero,K) : K \in \CA\}$ where $\CA$  is a compact set in $\CL_1^{+,s}$. 
If  \ref{H2} holds with $r<2$,  
then for every $R>0$ and $\la>0$ there is a compact 
$\CA_R \subset \CL_1^{+,s}$ such that, for all $N$, 
\[\sup_{\gamma \in \mathscr{A}}
 \int_{H^N} \II  \Big  ( \frac1{\la N} \sum_{i=1}^N C_{h_i}  \not \in \CA_R  \Big ) \gamma^{\otimes N} (dh_1 \dots dh_N)    \le e^{-RN}.
 \]
\end{lemma}

\begin{proof}
   Let  $\mathscr{A}^*=\{ \gamma \#  {C}_\cdot: \gamma \in \mathscr{A}\}$. 
 By Lemma \ref{Lemma:compactness2} for every $t$
  \[
\sup_{\mu \in \mathscr{A}^*}   \int_{\CL_1}   e^{t \| e  \|_1} \mu(de) =\sup_{\gamma \in \mathscr{A}}   \int_H  e^{t \| C_h  \|_1} \gamma(dh) <+\infty
\]
and  by  Lemma \ref{Lemma:compactness} $\mathscr{A}$ is tight. Since $f \mapsto C_f$ is continuous, 
the continuous mapping theorem and Prohorov theorem yields that 
also  $\mathscr{A}^*$ is tight. 
	Then, Proposition \ref{thm:deAcosta} yields some compact, convex, well balanced $V \subset \CL_1$ such that
\begin{equation}\label{boundA_1}
\sup_{\gamma \in \mathscr{A}}   \int_H  e^{ q_V( C_h )} \gamma(dh) =\sup_{\mu \in \mathscr{A}^*} \int_{\CL_1} e^{ q_V (e) } \mu(de) 
=:M<+\infty.
\end{equation}
Now, fix  $t_M>0$ and note that   
 $\{x \not \in a V\}=\{x : q_V(x) >a\}$.  Then,
\[
\begin{split}
& \int_{H^N} \II  \Big  ( \frac1{\la N} \sum_{i=1}^N C_{h_i}  \not \in t_M V   \Big ) \gamma^{\otimes N} (dh_1 \dots dh_N)  \\
&= \int_{H^N} \II  \Big  (  \sum_{i=1}^N C_{h_i} \not  \in \la  N t_M V   \Big ) \gamma^{\otimes N} (dh_1 \dots dh_N)  \\
&= \int_{H^N} \II  \Big  (  q_V (\sum_{i=1}^N C_{h_i} )  > \la N t_M \Big ) \gamma^{\otimes N} (dh_1 \dots dh_N)  \\
&\le  \int_{H^N}  e^{- \la N t_M} e^{q_V (\sum_{i=1}^N C_{h_i} )  } \gamma^{\otimes N} (dh_1 \dots dh_N)  \\
\end{split}
\]
where in the last step we use Chebyshev inequality. Then, $q_V(\sum_{i=1}^N C_{h_i}) \le   \sum_{i=1}^N  q_V(C_{h_i})$
by  subadditivity of $q_V$.
Hence, taking $t_M=\la^{-1}(R+\log(M))$,  we set 
$\CA_R\colon =t _M V \cap \CL_1^{+,s}$ and using also  \eqref{boundA_1} one gets 
\begin{equation}\label{bound_unifA1}
\begin{split}
\sup_{\gamma \in \mathscr{A}}
 \int_{H^N} \II  \Big  ( \frac1{\la N} \sum_{i=1}^N C_{h_i}  \not \in t_M V   \Big ) \gamma^{\otimes N} (dh_1 \dots dh_N) &   \le 
 \sup_{\gamma \in \mathscr{A}} e^{-t_M \la N}  \Big ( \int_H   e^{q_V (C_h )  } \gamma(dh) \Big)^N  \\
&  \le  e^{-t_M \la N +\log(M)} = e^{-RN}. \\
 \end{split} 
\end{equation}
\end{proof}

We are ready now to prove Proposition~\ref{Prop:exp.tight}.
 
 \begin{proof}[Proof of Proposition~\ref{Prop:exp.tight}]
For each fixed $R>0$, we will choose a special form for the corresponding compact set: $\mathcal{B}_R=\CA^{(2)}_R\times\dots \times \CA^{(L+1)}_R$, where for $\ell=2,\dots,L+1$, each $\CA^{(\ell)}_R$ is a compact subset of $\CLL$.
 The first step is to prove that $(\bK^2_{N_1})_N$ is exponentially tight, i.e. to find a compact set $\CA^{(2)}_R\subset \CLL$ such that, for $N$ sufficiently large,
 \[
 \PP\big(  \bK^2_{N_1}   \notin \CA^{(2)}_R\big)\leq \ee^{-NR}.
 \] 
 By \eqref{K2-sum}
and Lemma \ref{Lemma:compactness3} applied to 
 $\CA^{(1)}=\{\bK^1_{N_0}\}$ and $\la=\la_1$, 
 one gets that for every $R>0$ there exists a compact
 set $\CA_R=\CA^{(2)}_R$ such that 
 \[
 \PP\big (  \bK^2_{N_1}   \notin \CA^{(2)}_R\big)
= \int_{H^{N_1}} \II  \Big   ( \frac1{\la_1N_1} \sum_{i=1}^{N_1} C_{h_i}  \notin \CA^{(2)}_R   \Big ) \gamma_{\bK^1_{N_0}}^{\otimes N_1} (dh_1 \dots dh_{N_1}) \le e^{-RN_1}.
 \]
Now, one applies again Lemma \ref{Lemma:compactness3} to 
$\CA=\CA^{(2)}_R$, that is for  
$\mathscr{A}=\mathscr{A}_2  = \{ \CN_H(\bzero,K) : K \in \CA^{(2)}_R\}$,  to get 
that there is a compact set $\CA^{(3)}_R$ such that 
\begin{equation}\label{bound_unif_A1}
\sup_{\gamma \in \mathscr{A}_2}
 \int_{H^{N_2}} \II  \Big  ( \frac1{\la_1 N_2} \sum_{i=1}^{N_2} C_{h_i}  \not \in \CA^{(3)}_R\Big ) \gamma^{\otimes N_2} (dh_1 \dots dh_{N_2})    \le e^{-RN_2}.
\end{equation}
Now, note that 
\[
\PP\big((\bK^2_{N_1}, \bK^3_{N_2})  \not \in \CA^{(2)}_R \times \CA^{(3)}_R\big)
\le \PP\big( \bK^2_{N_1} \not \in \CA^{(2)}_R\big) +\PP\big(  \bK^2_{N_1}  \in \CA^{(2)}_R,  \bK^3_{N_2}   \not \in  \CA^{(3)}_R\big).
\]
Using \eqref{transitionK}, one can write
\[
\begin{split}
\PP\big( \bK^2_{N_1} &  \in \CA^{(2)}_R,  \bK^3_{N_2}   \not \in  \CA^{(3)}_R\big)
=\E\big[  \II( \bK^2_{N_1}  \in \CA^{(2)}_R ) \E[  \II (  \bK^3_{N_2}   \not \in  \CA^{(3)}_R) | \bK^2_{N_1} ] \big] \\
& =\E\Big [  \II( \bK^2_{N_1}  \in \CA^{(2)}_R )
 \int_{H^{N_2}} \II  \Big   ( \frac1{\la_2N_2} \sum_{i=1}^{N_2} C_{h_i}  \not \in t_M V   \Big ) \gamma_{\bK^2_{N_1}}^{\otimes N_2} (dh_1 \dots dh_{N_2}) \Big ]
\\
\end{split}
\]
where $\gamma_K\colon =\CN_H(\bzero,K) $. By \eqref{bound_unifA1} and the definition of $\mathscr{A}_2$, this becomes 
\[
\PP\big(  \bK^2_{N_1}   \in {\CA}^{(2)}_R,  \bK^3_{N_2}   \not \in  \CA^{(3)}_R\big)
\le 
\sup_{\gamma \in \mathscr{A}_1}
 \int_{H^{N_2}} \II  \Big  ( \frac1{\la_2 {N_2}} \sum_{i=1}^{N_2} C_{h_i}  \not \in t_M V   \Big ) \gamma^{\otimes N_2} (dh_1 \dots dh_{N_2})    \le e^{-RN_2} .
  \]
Combining all together   one gets
\[
\limsup_{N\to\infty}\frac 1N \log \PP\big((\bK^2_{N_1}, \bK^3_{N_2})  \not \in \CA^{(2)}_R \times \CA^{(3)}_R\big) \le -R', 
\]
where $R'=R\max\{m_1,m_2\}$. Iterating this procedure a finite number of times, one obtains  the thesis. 
 \end{proof}

%
 %SS EXPTIGHT C0
 %
\subsection{Exponential tightness of the law  in  $\CC^{+,s}\times\dots\times \CC^{+,s}$: proof of Proposition~\ref{prop:exp_tight_C0}}
\label{ss:expTight_C_proof}

Recall that for the proof of Proposition~\ref{prop:exp_tight_C0}, we assume additionally~\ref{H4}. To prove the result, we shall consider 
a compact set $\CC_M=\CC^{(2)}_M\times\dots\times \CC^{(L+1)}_M$, where for any $\ell\leq L+1$ the set $\CC^{(\ell)}_M\subset \CC^{+,s}$ is constructed as follows: 
there exist $C_\ell, C_\ell^{\prime}<\infty$ such that $\CK^{\ell}_{N_{\ell-1}}\in \CC^{(\ell)}_M$ if and only if
\begin{align*}
    &\sup_{(\bx,\by),(\bx',\by')\in\UU^2}\frac{(\CK^{\ell}_{N_{\ell-1}}(\bx,\by)-\CK^{\ell}_{N_{\ell-1}}(\bx',\by'))^2}{\|\bx-\bx'\|^2+\|\by-\by'\|^2}\le C_\ell;\\
    &\\
    &\sup_{(\bx,\by)\in\UU^2}(\CK^{\ell}_{N_{\ell-1}}(\bx,\by))^2\le C_\ell^{\prime}.
\end{align*}
As already specified earlier, $\CK^{\ell}_{N_{\ell-1}}\in \CC^{+,s}$ a.s..

We will prove Proposition~\ref{prop:exp_tight_C0} by means of two lemmas, for which we need to define some events first. For every $\ell=1,\dots,L$, let $A^{(\ell)}<\infty$ and 
\begin{equation}
    \label{eq:equicont}
\mathcal{A}_{(\ell)}\colon =\Big\{\sup_{\bx,\by\in \UU} \frac{\sum_{i=1}^{N_\ell}(h^{(\ell)}_i(\bx)-h^{(\ell)}_i(\by))^2}{\|\bx-\by\|^2}\le  A^{(\ell)}N_\ell\Big\}.
\end{equation}
Similarly, for $B^{(\ell)}<\infty$, we define the event
\begin{equation}
    \label{eq:boundedness}
\mathcal B_{(\ell)}\colon=\Big\{\sup_{\bx\in\UU}\sum_{j=1}^{N_\ell}|h^{(\ell)}_j(\bx)|^2\le B^{(\ell)}N_\ell\Big\}.
\end{equation}

\begin{lemma}\label{C0tightLemma1} Assume that $\UU \subset \RE^{N_0}$ is compact and  \ref{H1}-\ref{H2}-\ref{H3}-\ref{H4}, with
$r \leq 2$ in \ref{H2}. Fix $\ell=1,\dots, L$, then for any $M<\infty$, there exist $A^{(\ell)}, B^{(\ell)} < \infty$ such that
\[
\limsup_{N\to\infty}\frac 1{N_\ell}\log \PP\left(\mathcal{A}_{(\ell)}^{\mathsf{C}}\cup\mathcal{B}_{(\ell)}^{\mathsf{C}}\right)
\le -M,
\]
where $\mathcal{A}_{(\ell)}$ and $\mathcal{B}_{(\ell)}$ are the events defined in \eqref{eq:equicont} and \eqref{eq:boundedness}, respectively.
\end{lemma}

\begin{lemma}\label{C0tightLemma2} Assume that $\UU \subset \RE^{N_0}$ is compact and  \ref{H1}-\ref{H2}-\ref{H3}-\ref{H4}, with
$r\leq 2$ in \ref{H2}. Fix $\ell=2,\dots, L+1$, then under the event $\mathcal{A}_{(\ell-1)}\cap\mathcal{B}_{(\ell-1)}$ (defined in \eqref{eq:equicont} and \eqref{eq:boundedness}) there exist $C_\ell, C_\ell^{\prime}<\infty$ such that 
\begin{align*}
    &\sup_{(\bx,\by),(\bx',\by')\in\UU^2}\frac{(\CK^{\ell}_{N_{\ell-1}}(\bx,\by)-\CK^{\ell}_{N_{\ell-1}}(\bx',\by'))^2}{\|\bx-\bx'\|^2+\|\by-\by'\|^2}\le C_\ell;\\
    &\\
    &\sup_{(\bx,\by)\in\UU^2}\CK^{\ell}_{N_{\ell-1}}(\bx,\by)^2\le C_\ell^{\prime}.
\end{align*}
    
\end{lemma}

\begin{proof}[Proof of Proposition 
\ref{prop:exp_tight_C0}]
Combine Lemma~\ref{C0tightLemma1} and Lemma~\ref{C0tightLemma2}. 
\end{proof}

Now, we prove Lemma~\ref{C0tightLemma1} and Lemma~\ref{C0tightLemma2}.

\begin{proof}[Proof of Lemma \ref{C0tightLemma1}]
	We prove the claim by induction on $\ell$.

	\emph{Induction start, $\ell=1$.} Fix $\bx,\by\in\UU$, 
because of the easy form of $h^{(1)}(\mathbf{x})$, 
\[
\begin{split}
\frac1{N_{1}}\sum_{i=1}^{N_{1}}(h^{(1)}_i(\bx)-h^{(1)}_i(\by))^2
&= 
\frac1{N_1}\sum_{i=1}^{N_1}\frac1{N_0}\Big (\sum_{j=1}^{N_0}W_{ij}^{(0)}(x_j-y_j)\Big)^2 \\
& \le \frac{1}{N_1}\sum_{i=1}^{N_1}\left(\frac 1{N_0}\sum_{j=1}^{N_0}(W_{ij}^{(0)})^2\right)\left(\sum_{j=1}^{N_0}(x_j-y_j)^2\right), \\
\end{split}
\]
where the inequality is due to Cauchy-Schwarz. Then, we have that
\begin{align*}
\PP\left(\sup_{x,y\in \UU} \frac{\sum_{i=1}^{N_1}\frac1{N_0}(\sum_{j=1}^{N_0}W_{ij}^{(0)}(x_j-y_j))^2}{\|\bx-\by\|^2}\ge N_1 A^{(1)}\right)&\le \PP\left(\sum_{i=1}^{N_1}\left(\frac 1{N_0}\sum_{j=1}^{N_0}(W_{ij}^{(0)})^2\right)\ge N_1 A^{(1)}\right)\\
&\le \ee^{-N_1A^{(1)} t}\left(\mathbb{E}[\ee^{t\mathbf{X^{(0)}}}]\right)^{N_1},
\end{align*}
where  
$\mathbf{X^{(0)}}\overset{d}{=}\frac 1{N_0}\sum_{j=1}^{N_0}(W_{1j}^{(0)})^2$ has distribution $\Gamma(N_0/2,\la_0N_0/2)$, since $W_{ij}^{(0)}$ are i.i.d. random variables with law $\CN(0,\la_0^{-1})$, which in particular has finite exponential moment $\mathbb{E}[\ee^{t\mathbf{X^{(0)}}}]$ if $t$ is small enough. Notice that above we used that $(\frac 1{N_0}\sum_{j=1}^{N_0}(W_{ij}^{(0)})^2)_{i=1,\dots,N_1}$ are i.i.d. random variable (with the same distribution as $\mathbf{X^{(0)}}$).  One can then choose $A^{(1)}$ large enough to have 
\[\ee^{-t A^{(1)}}\mathbb{E}[\ee^{t\mathbf{X^{(0)}}}]\le \ee^{-M}.
\]
This shows that 
$
\PP\big(\mathcal{A}_{(1)}^{\mathsf{C}} \big)
\le  e^{-MN_1}.
$
Similarly, 
\[
\PP\Big(\sup_{x\in \UU} \sum_{i=1}^{N_1}(\sum_{j=1}^{N_0}W_{i,j}^{(0)}x_j)^2\ge B^{(1)}N_0N_1\Big)\le \ee^{-N_1\frac t{U}B^{(1)}}\Big(\mathbb{E}[\ee^{t\mathbf{X^{(0)}}}]\Big)^{N_1},
\]
with $U\colon=\sup_{\mathbf{x}\in\UU}\sum_{j=1}^{N_0}(x_j)^2=\sup_{\mathbf{x}\in\UU}\|\bx\|^2$. Hence, one can  choose $B^{(1)}$ large enough such that 
\[\ee^{-\frac tUB^{(1)}}\mathbb{E}[\ee^{t\mathbf{X^{(0)}}}]\le \ee^{-M}.\]
Which give 
$\P(\mathcal{B}_{(\ell)}^{\mathsf{C}}) \leq \ee^{-N_1 M}.$
This proves the case $\ell=1$. 

%%%%%%%%%%%%%%%%%%%%%%%%%%%%%%%%%
%%%%%%%%%%%%%%%%%%%%%%%%
\emph{Induction step.}
Now, suppose that we proved the statement for $\ell-1$, let us prove it for $\ell$. 
To simplify the notations 
we write $W_{ij}$ in place of $W_{ij}^{(\ell-1)}$.  
Using  the explicit expression of $h^{(\ell)}$ in terms of $h^{(\ell-1)}$, see~\eqref{main_recursion}, 
we write
\begin{align*}
    \frac1{N_\ell}& \sum_{i=1}^{N_\ell}(h^{(\ell)}_i(\bx)-h^{(\ell)}_i(\by))^2=\frac1{N_\ell}\sum_{i=1}^{N_\ell}{\frac1{N_{\ell-1}}}\left(\sum_{j=1}^{N_{\ell-1}}W_{ij}\left(\s(h^{(\ell-1)}_j(x))-\s(h^{(\ell-1)}_j(y))\right)\right)^2\\
    =&\frac1{N_{\ell-1}}\sum_{j,j'}\left(\s(h^{(\ell-1)}_j(x))-\s(h^{(\ell-1)}_j(y))\right)\left(\s(h^{(\ell-1)}_{j'}(x))-\s(h^{(\ell-1)}_{j'}(y))\right)\frac1{N_\ell}\sum_{i=1}^{N_\ell}W_{ij}W_{ij'}\\
    =& \frac1{N_{\ell-1}}\langle \mathbf{\s^h(x,y)}, \mathbf{U}\mathbf{\s^h(x,y)}\rangle
\end{align*}
where we indicate with $\mathbf{U}$ the $N_{\ell-1}\times N_{\ell-1}$ matrix with entries as follows:
\[
U_{jj'}\colon = \frac1{N_\ell}\sum_{i=1}^{N_\ell}W_{ij}W_{ij'}
\]
and with $\mathbf{\s^h(x,y)}$ the $N_{\ell-1}$-dimensional vector \[\mathbf{\s^h(\bx,\by)}=(\s(h_1^{(\ell-1)}(\bx))-\s(h_1^{(\ell-1)}(\by)),\dots, \s(h_{N_{\ell-1}}^{(\ell-1)}(\bx))-\s(h_{N_{\ell-1}}^{(\ell-1)}(\by))^\top.\]
Notice that $\mathbf{U}$ is a symmetric matrix with real entries, that can be expressed as $\frac1{{N_\ell}}\mathbf{W}^\top\mathbf{W}$ where $\mathbf{W}$ is the $N_\ell\times N_{\ell-1}$ matrix with entries $W_{ij} \stackrel{iid}{\sim} \CN(0,\la_{\ell-1}^{-1})$.

In order to bound the quantity $\frac1{N_{\ell-1}}\langle \mathbf{\s^h(x,y)}, \mathbf{U}\mathbf{\s^h(x,y)}\rangle$ we can use Rayleigh quotient to write 
\[
\langle \mathbf{\s^h(x,y)}, \mathbf{U}\mathbf{\s^h(x,y)}\rangle\le \la(\mathbf{U})\langle \mathbf{\s^h(x,y)}, \mathbf{\s^h(x,y)}\rangle, 
\]
where $\la(\mathbf{U})$ is the largest eigenvalue of the matrix $\mathbf{U}$. 
Setting
\[
\Sigma^{(\ell)}(\bx,\by)\colon=\frac{1}{N_{\ell-1}\|\bx-\by\|^2}\sum_{j=1}^{N_{\ell-1}}\left(\s(h^{(\ell-1)}_j(\bx))-\s(h^{(\ell-1)}_j(\by))\right)^2.
\]
one has 
\[
\sup_{\bx,\by \in \UU}
\frac{
\sum_{i=1}^{N_\ell}(h^{(\ell)}_i(\bx)-h^{(\ell)}_i(\by))^2}{N_\ell\|\bx-\by\|^2}\le
\la(\mathbf{U})
\sup_{\bx,\by \in \UU}
\frac{\|\mathbf{\s^h}(\bx,\by)\|^2}{N_{\ell-1}\|\bx-\by\|^2}
=\la(\mathbf{U})\sup_{\bx,\by \in \UU}\Sigma^{(\ell)}(\bx,\by).
\]
Thanks to \ref{H4}, we see that
\[
\sum_{j=1}^{N_{\ell-1}}\left(\s(h^{(\ell-1)}_j(\bx))-\s(h^{(\ell-1)}_j(\by))\right)^2\leq L_{\s^2}\sum_{j=1}^{N_{\ell-1}}\left(h^{(\ell-1)}_j(\bx)-h^{(\ell-1)}_j(\by)\right)^2,
\]
hence on the event $\mathcal{A}_{(\ell-1)}$, we have that
\[
\sup_{h^{(\ell-1)}\in \mathcal{A}_{(\ell-1)}}\sup_{\bx,\by\in\UU}\Sigma^{(\ell)}(\bx,\by)\le \Sigma^{(\ell)}\colon =L_\s^2 A^{(\ell-1)}<\infty.
\]

Being $\mathbf{U}=\frac1{{N_\ell}}\mathbf{W}^\top\mathbf{W}$, its largest eigenvalue is equal to the square of the largest \textit{singular value} of the matrix $\frac 1{\sqrt{N_\ell}}\mathbf{W}$.
Hence, we can use \cite[Theorem 4.4.5 ]{Ver19} (see Proposition \ref{thmVer19}
in the Appendix) to obtain the
following inequality: for every $t>0$  
\begin{equation}\label{bound-sing-val}
\PP\Big({\la(\mathbf{U})}\ge C\Big(1+\sqrt{{N_{\ell-1}}/{N_\ell}}+t\Big)^2\Big)\le 2\ee^{-N_\ell t^2}
\end{equation}
where $C=C_{\ell-1}$ is a given constant (independent of $N_{\ell-1}$ and $N_\ell$).

Combining the previous results, given $M>0$, one can find $A^{(\ell)}$ such that 
\[
\PP\Big(\sup_{\bx,\by\in \UU} \frac{\sum_{i=1}^{N_\ell}(h^{(\ell)}_i(x)-h^{(\ell)}_i(y))^2}{N_\ell\|\bx-\by\|^2}\ge A^{(\ell)}, \mathcal{A}_{(\ell-1)}\Big)
\le \PP( \la(\mathbf{U}^{\ell-1}) \ge A^{(\ell)}/\Sigma^{(\ell)} )
\le 2 \ee^{-M N_\ell}.
\]
In this way we have 
\begin{align*}   
\PP&\Big(\sup_{\bx,\by\in \UU} \frac{\sum_{i=1}^{N_\ell}(h^{(\ell)}_i(\bx)-h^{(\ell)}_i(\by))^2}{N_\ell\|\bx-\by\|^2}\ge A^{(\ell)}\Big)\\
&\qquad\le \PP\Big(\sup_{x,y\in \UU} \frac{\sum_{i=1}^{N_\ell}(h^{(\ell)}_i(\bx)-h^{(\ell)}_i(\by))^2}{{N_\ell}\|\bx-\by\|^2}\ge A^{(\ell)}, \mathcal{A}_{(\ell-1)}\Big)+\PP(\mathcal{A}_{(\ell-1)}^{\mathsf{C}})\le 2\ee^{-MN_\ell}+\PP(\mathcal{A}_{(\ell-1)}^{\mathsf{C}})
\end{align*}
which  gives
$
\limsup_{N\to\infty}\frac 1{N_\ell}\log \PP\big(\mathcal{A}_{(\ell)}^{\mathsf{C}} \big)
\le -M.$
The bound for $\mathcal{B}_{(\ell)}$ is obtained similarly. 
Recalling that $W_{ij}=W_{ij}^{(\ell-1)}$, 
one can write 
\begin{align*}
    \frac1{N_\ell}\sum_{i=1}^{N_\ell}(h_i^{(\ell)}(\bx))^2&\overset{\eqref{main_recursion}}{=} \frac1{N_\ell}\sum_{i=1}^{N_\ell}\frac1{N_{\ell-1}}\Big(\sum_{j=1}^{N_{\ell-1}}W_{ij}\s(h^{(\ell-1)}_j(\bx))\Big)^2\\
    &= \frac1{N_{\ell-1}}
    \sum_{j=1}^{N_{\ell-1}}\sum_{j'=1}^{N_{\ell-1}}
    \s(h^{(\ell-1)}_j(\bx)) \s(h^{(\ell-1)}_{j'}(\bx))
    \frac1{N_\ell} \sum_{i=1}^{N_\ell}  W_{ij}W_{ij'} \\
    & \le \la(\mathbf{U}) 
    \frac1{N_{\ell-1}} \sum_{j=1}^{N_{\ell-1}} \s(h^{(\ell-1)}_j(\bx))^2.
\end{align*}
On $\mathcal{B}_{(\ell-1)}$, we have
\begin{align*}
    \sup_{\bx\in\UU}\frac1{N_{\ell-1}}\sum_{i=1}^{N_{\ell-1}}\s(h_i^{(\ell-1)}(\bx))^2
    &\overset{\ref{H2}}{\le} \sup_{\bx\in\UU} \Big(A+A\frac1{N_{\ell-1}}\sum_{j=1}^{N_{\ell-1}}|h^{(\ell-1)}_j(\bx)|^{r}\Big)\\
    &\overset{Jensen}{\le} \sup_{\bx\in\UU}\Big(A+A\Big(\frac1{N_{\ell-1}}\sum_{j=1}^{N_{\ell-1}}|h^{(\ell-1)}_j(\bx)|^2\Big)^{r/2}\Big)\\
	&\overset{\mathcal{B}_{(\ell-1)}}{\le} A(1+(B^{(\ell-1)})^{r/2})<\infty,
\end{align*}
where we use the growth condition on $\s$ and Jensen inequality, since $r/2 \leq 1$. 
Applying once again \eqref{bound-sing-val}, one can find a constant $B^{(\ell)}$ such that 
$
\PP \big(  \mathcal{B}_{(\ell)}^{\mathsf{C}}, \mathcal{B}_{(\ell-1)}\big)
\leq 2\ee^{-N_\ell M}.
$
As before this gives 
$
\limsup_{N\to\infty}\frac 1{N_\ell}\log \PP\big(\mathcal{B}_{(\ell)}^{\mathsf{C}} \big)
\le -M.
$
\end{proof}

\begin{proof}[Proof of Lemma \ref{C0tightLemma2}]
    We use the explicit formulation in \eqref{eq:K_in_C0}. For any fixed $\ell=2,\dots, L+1$, we have the following
    \begin{equation}\label{eq:CK_unif_cont}
        \begin{split}
        (\CK^{\ell}_{N_{\ell-1}}(\bx,\by)-\CK^{\ell}_{N_{\ell-1}}(\bx',\by'))^2&=\frac1{N^2_{\ell-1}}\Big(\sum_{i=1}^{N_{\ell-1}}\s(h^{(\ell-1)}_i(\bx))\s(h^{(\ell-1)}_i(\by))-\s(h^{(\ell-1)}_i(\bx'))\s(h^{(\ell-1)}_i(\by'))\Big)^2\\
        &\le\frac2{N^2_{\ell-1}}\Big(\sum_{i=1}^{N_{\ell-1}}\s(h^{(\ell-1)}_i(\bx))\Big(\s(h^{(\ell-1)}_i(\by))-\s(h^{(\ell-1)}_i(\by'))\Big)\Big)^2\\
        &\qquad+\frac2{N^2_{\ell-1}}\Big(\sum_{i=1}^{N_{\ell-1}}\s(h^{(\ell-1)}_i(\by'))\Big(\s(h^{(\ell-1)}_i(\bx))-\s(h^{(\ell-1)}_i(\bx'))\Big)\Big)^2.
    \end{split}
    \end{equation}
    By symmetry, we focus on the first term of the above sum, which we bound as follows, using Cauchy-Schwarz inequality:
    \begin{align*}
        \frac1{N^2_{\ell-1}}&\Big(\sum_{i=1}^{N_{\ell-1}}\s(h^{(\ell-1)}_i(\bx))\Big(\s(h^{(\ell-1)}_i(\by))-\s(h^{(\ell-1)}_i(\by'))\Big)\Big)^2\\
        &\overset{CS}{\le} \Big(\frac1{N_{\ell-1}}\sum_{i=1}^{N_{\ell-1}}\s(h^{(\ell-1)}_i(\bx))^2\Big)\Big(\frac1{N_{\ell-1}}\sum_{i=1}^{N_{\ell-1}}\Big(\s(h^{(\ell-1)}_i(\by))-\s(h^{(\ell-1)}_i(\by'))\Big)^2\Big).
    \end{align*}
 On the event  $\mathcal{B}_{(\ell-1)}$ one has 
\begin{align*}
    \frac1{N_{\ell-1}}\sum_{i=1}^{N_{\ell-1}}\s(h^{(\ell-1)}_i(\bx))^2&\overset{\ref{H2}}{\le} \frac1{N_{\ell-1}}\sum_{i=1}^{N_{\ell-1}}A(1+(h^{(\ell-1)}_i(x))^r)\\
    &\le A+A \frac1{N_{\ell-1}}\sum_{i=1}^{N_{\ell-1}}(h^{(\ell-1)}_i(\bx))^{r}\\
    &\overset{Jensen}{\le} A+A\Big(\frac1{N_{\ell-1}}\sum_{i=1}^{N_{\ell-1}}(h^{(\ell-1)}_i(\bx))^2\Big)^{r/2}\overset{\mathcal{B}_{(\ell-1)}}{\le} A(1+(B^{(\ell-1)})^{r/2}),
\end{align*}
where we used assumption (H2), Jensen inequality (since $r/2 \leq 1$) and the definition of the event $\mathcal{B}_{(\ell-1)}$. To bound the second term, we proceed as follows:
\begin{align*}
    \frac1{N_{\ell-1}}\sum_{i=1}^{N_{\ell-1}}\Big(\s(h^{(\ell-1)}_i(\by))-\s(h^{(\ell-1)}_i(\by'))\Big)^2&\overset{(H4)}{\le} L^2_{\s} \frac1{N_{\ell-1}}\sum_{i=1}^{N_{\ell-1}}\Big(h^{(\ell-1)}_i(\by)-h^{(\ell-1)}_i(\by')\Big)^2\\
    &\le L^2_{\s} A_{(\ell-1)}\|\by-\by'\|^2,
\end{align*}
where we used assumption (H4) and the fact that we are on the event $\mathcal{A}_{(\ell-1)}$. Inserting the two upper bounds above in \eqref{eq:CK_unif_cont}, we see that 
\begin{align*}
    \sup_{(\bx,\by),(\bx',\by')\in\UU^2}\frac{(\CK^{\ell}_{N_{\ell-1}}(\bx,\by)-\CK^{\ell}_{N_{\ell-1}}(\bx',\by'))^2}{\|\bx-\bx'\|^2+\|\by-\by'\|^2} \overset{\mathcal{A}_{(\ell-1)}}{\le}  4A(1+(B^{(\ell-1)})^{r/2})  L^2_{\s} A_{(\ell-1)}=\colon C_\ell<\infty.
\end{align*}
We repeat some of the arguments above to see that
\begin{align*}
  N^2_{\ell-1}  \CK^{\ell}_{N_{\ell-1}}(\bx,\by)^2&\overset{\eqref{eq:K_in_C0}}{=}\Big(\sum_{i=1}^{N_{\ell-1}}\s(h^{(\ell-1)}_i(\bx))\s(h^{(\ell-1)}_i(\by))\Big)^2\\
    &\overset{CS}{\le} \Big(\sum_{i=1}^{N_{\ell-1}}\s(h^{(\ell-1)}_i(\bx))^2\Big)\Big(\sum_{i=1}^{N_{\ell-1}}\s(h^{(\ell-1)}_i(\by))^2\Big)\\
    &\overset{\ref{H2}}{\le} \Big(\sum_{i=1}^{N_{\ell-1}}A(1+(h^{(\ell-1)}_i(\bx))^r)\Big)\Big(\sum_{i=1}^{N_{\ell-1}}A(1+(h^{(\ell-1)}_i(\by))^r)\Big)\\
    &\overset{\substack{Jensen,\\ \mathcal{B}_{(\ell-1)}}}\le   N^2_{\ell-1}\Big(A(1+(B^{(\ell-1)})^{r/2})\Big)^2.
\end{align*}
	Hence, dividing both sides by $N^2_{\ell -1}$ concludes the proof.
\end{proof}

\subsection{Proof of Proposition~\ref{prop:posterior}}\label{sec:proof_posterior}
Let
	$\bS_{1:P} =\vvec[h^{(L+1)}(\bx_1),\dots,h^{(L+1)}(\bx_P)]$. 
By \eqref{posterior_infty}
the posterior of $\boldsymbol{\CK}_N$
can be written as  
\[
\PP(\boldsymbol{\CK}_N \in A| \by_{1:P})
=\frac{\int_{A \times \RE^{PD}}  \CL(\by_1,\dots,\by_P| \bs_{1},\dots,\bs_P) P_{N,\mathrm{prior}}(d \overline{\CKK} d\bs_{1:P})}{ \int \CL(\by_1,\dots,\by_P|\bs_{1},\dots,\bs_P ) P_{N,\mathrm{prior}}(d\overline{\CKK} d\bs_{1:P})}
\qquad A \in \CBB\big((\CC^{+,s})^L\big)
\]
where $\overline{\CKK}=(\CKK^2,\dots,\CKK^{L+1})$, $\bs_{1:P}=\vvec[\bs_1,\dots,\bs_{P}]$
and 
$P_{N,\mathrm{prior}}(d\overline{\CKK}  ds_{1:P})$ is the distribution 
of $(\boldsymbol{\CK}_N,\bS_{1:P})$ (induced by $P_{N,\mathrm{prior}}(d\theta)$).
Hence, in order to compute the posterior of 
$\boldsymbol{\CK}_N$, it is enough to consider
the joint law of $(\boldsymbol{\CK}_N,\bS_{1:P},\bY_{1:P})$. 
Now, note that 
\begin{equation}\label{condSgivenK}
   \bS_{1:P} | \CK^{L+1}_N 
\sim \CN(\bzero,\Sigma)  
\end{equation}
	where $\Sigma: = \Sigma(\CK^{L+1}_N)$ has been defined in \eqref{eq:sig}. To go further, we expand the joint law 
of $(\boldsymbol{\CK}_N,\bS_{1:P},\bY_{1:P})$
with an auxiliary variable $\mathbf{Z}_{1:P}$ which simplifies computations. 
To this end, let us consider the  joint distribution in $(\by_{1:P},\bs_{1:P},\bz_{1:P},\overline{\CKK})$ given by 
\begin{equation}\label{joint}
\begin{split}
\mu(d\overline{\CKK} d\bz_{1:P}& d\bs_{1:P} d\by_{1:P})  =
\mu(d\by_{1:p}|\bz_{1:P},\overline{\CKK}) \mu(d\bs_{1:p}|\bz_{1:P},\overline{\CKK}) \mu(d\bz_{1:P} )\mu(d\overline{\CKK})  \\
& 
:=\delta_{\Sigma^{1/2} \mathbf{z}_{1:P}}(d\bs_{1:P})\CN(d\by_{1:p}| \Sigma^{1/2} \bz_{1:P},\beta^{-1} \eI_{DP} ) 
\CN(d\bz_{1:P}| \mathbf{0},   \eI_{DP} )  \mathcal{Q}_N(d\overline{\CKK}),
\end{split}
\end{equation}
where now $\Sigma=\Sigma(\CKK^{L+1})$.
Note that above $\mathcal{Q}_N(d\overline{\CKK})=\mu(d\overline{\CKK})$, $\bz_{1:P}$ and $\by_{1:P}$ are in $\RE^{DP}$, $\bs_{1:P}=\Sigma^{1/2} \mathbf{z}_{1:P}$ and 
$\mathbf{z}_{1:P}$ and $\overline{\CKK}$ are independent. 
Hence,
 $\mu(d\bs_{1:P}|\overline{\CKK})$
 is a Gaussian distribution
with mean  $\mathbf{0}$ and covariance matrix  $\Sigma$ and by \eqref{condSgivenK} it follows that 
 $\mu(d\overline{\CKK} d\bs_{1:p} d\by_{1:P})$ is the 
joint law of $(\boldsymbol{\CK}_N,\bS_{1:P},\bY_{1:P})$.
In order to derive the conditional 
distribution of $\boldsymbol{\CK}_N$ 
given $\bY_{1:P}$ we first note that 
\begin{equation}\label{zQ|y}
\mu(d \bz_{1:P} d\overline{\CKK}| \by_{1:P} ) =  \mu(  d \bz_{1:P} |\overline{\CKK}, \by_{1:P})\mu(d\overline{\CKK}|  \by_{1:P}) \propto f(\bz_{1:P},\by_{1:P}| \overline{\CKK} )  \mathcal{Q}_N(d\overline{\CKK}) d\bz_{1:P} 
\end{equation}
with 
\[
f(\bz_{1:P},\by_{1:P}| Q ):=   
{e^{-\frac{\beta}2 (\Sigma^{1/2} \bz_{1:P} -\by_{1:P})^\top ( \Sigma^{1/2} \bz_{1:P}-\by_{1:P})} } 
e^{- \frac12\bz_{1:P}^\top \bz_{1:P} } .
\]
Setting  
\[
\mathbf{m}= \beta ( {\beta}\Sigma+\eI_{DP} )^{-1} \Sigma^{1/2} \by_{1:p}, 
\]
which is well-defined since $ {\beta}\Sigma+\eI_{DP}>0$, 
one checks that  
\[
\begin{split}
& \frac{\beta}2 (\Sigma^{1/2} \bz_{1:P}-\by_{1:P})^\top (\Sigma^{1/2} \bz_{1:P} -\by_{1:P})  - \frac12\bz_{1:P}^\top \bz_{1:P}  \\
&=\frac{\beta}2 \Big [  \by_{1:p} ^\top( \eI_{DP}- \beta  \Sigma^{1/2}  ({\beta}\Sigma+\eI_{DP} )^{-1}\Sigma^{1/2}  ) \by_{1:p}  \Big] \\
& +\frac12 (\bz_{1:P}-\mathbf{m})^\top ( {\beta} \Sigma+\eI_{DP})  (\bz_{1:P}-\mathbf{m})  .
\\
\end{split}
\]
Noticing  that 
$
( \eI_{DP}- \beta  \Sigma^{1/2}  ({\beta}\Sigma+\eI_{DP} )^{-1}\Sigma^{1/2}  )
= ({\beta}\Sigma+\eI_{DP} )^{-1},
$
one can write 
\[
f(\bz_{1:P},\by_{1:P}| \CKK )
={e^{-\frac12 \Psi (\CKK^{L+1}|\by_{1:P}) }} 
 \frac{ e^{-\frac12 (\bz_{1:P}-\mathbf{m})^\top ( {\beta} \Sigma+\eI_{DP})  (\bz_{1:P}-\mathbf{m})  }}{ \det(( {\beta} \Sigma+\eI_{DP})
 ^{-1})^{1/2} } 
\]
where
\[
\Psi (\CKK^{L+1}|\by_{1:P})=
\beta \by_{1:P}^\top ( \eI_{DP}+\beta \Sigma(\CKK^{L+1}))^{-1} \by_{1:P}
 +\log(\det(\eI_{DP} +\beta \Sigma(\CKK^{L+1})).
\]
Then  
\begin{equation}\label{cond_Z_givenQy}
\mu(d\bz_{1:P},d\overline{\CKK}|\by_{1:P})  =\frac{ e^{-\frac12 (\bz_{1:P}-\mathbf{m})^\top ( {\beta} \Sigma+\eI_{DP})  (\bz_{1:P}-\mathbf{m})  }}{(2 \pi)^{\frac{DP}2} \det(( {\beta} \Sigma+\eI_{DP})
 ^{-1})^{1/2} }  \mathcal{Q}_N(d\overline{\CKK}| \by_{1:P}) 
\end{equation}
with 
\[
 \mathcal{Q}_N(d\overline{\CKK} | \by_{1:P}) =% \mathcal{Q}(dQ |\tilde \bX, \lambda, \by_{1:P})   :=
 \frac{e^{-\frac12 \Psi (\CKK^{L+1}|\by_{1:P}) }\mathcal{Q}_N(d\overline{\CKK})  }{\int_{\CS^+_{D}}e^{-\frac12 \Psi (\CKK^{L+1}| \by_{1:P}) } \mathcal{Q}_N(d\overline{\CKK}) } . 
\]
Marginalizing with respect to $\bz_{1:P}$ gives $\mu(d\overline{\CKK}| \by_{1:P})=
\mathcal{Q}_N(d\overline{\CKK}| \by_{1:P})$ and hence the thesis.

\appendix

\section{Appendix}

\subsection{Hilbert-Schmidt and trace-class operators} \label{app:trace_class}
 Let $H$ be a separable  (real) Hilbert space
with scalar product $( \cdot, \cdot )_H$ 
with the usual identification $H'=H$. 
Denote by $\CL_\infty(H)$ the set of bounded linear operators with 
the operator norm 
$\|K\|=\sup_{f \in H:\|f\|_H=1} \|Kf\|_H$. 
A linear operator  $K : H \to H$ 
is positive (in symbols, $K \ge 0$)
if $\inf_{f \in H}( Kf, f )_H \ge 0$  and 
it is self-adjoint (symmetric) if $K^*=K$. 
If $K \ge 0$ then there is a unique positive 
operator $\sqrt{K}$ such that $\sqrt{K}\sqrt{K}=K$. 
Given 
$K$ one sets $|K|=\sqrt{K^*K}$
and $\tr(|K|)\colon =\sum_n (e_n,|K| e_n)_H$ where $(e_n)_n$ is any orthonormal basis, this does not depend on the choice of $(e_n)_n$. See VI.4 and VI.5 in \cite{ReedSimonVol1}. 
Denote by  $\CL_1(H)$ the set  of all  linear operators $K$ 
such that $\tr|K|<+\infty$. Such operators are  known as  trace-class operators or nuclear operators. 
If $K$
is in $\CL_1(H)$ 
then $\tr(K)\colon =\sum_n (e_n,K e_n)_H$ is absolutely convergent, where $(e_n)_n$ is any orthonormal basis. This does not depend on the basis,  \cite[Thm. VI.24]{ReedSimonVol1}. 
In point of fact   $\|K\|_1\colon =\tr|K|$  is a norm 
and one has the following.
\begin{enumerate}[label=(K\arabic*)]
\item \label{Fact0}
\textit{The space $(\CL_1(H),\|\cdot\|_1)$ is a Banach space.
If $H$
is separable, also $(\CL_1(H),\|\cdot\|_1)$ is separable. }
\end{enumerate}
See, e.g.,   \cite[Thm. VI.20]{ReedSimonVol1} and 
\cite[Thm. 18.11 (d)]{Conway}.

The operators such that $\|K\|_2^2\colon =\tr(K^*K)<+\infty$ are known as 
 Hilbert-Schmidt operators, and the set of such operators is 
  denoted by $\CL_2(H)$. One has $\CL_1(H) \subset \CL_2(H) \subset \CL_\infty(H)$ and $\CL_2(H)$ is contained in the space of compact operators.
 An operator  $K$ is in $\CL_1(H)$ if and only if 
$K=A^*B$ with $A$ and $B$ in $\CL_2(H)$.  
The space $(\CL_2(H),\|K\|_2)$ is an Hilbert space 
with scalar product $(A,B)_2=\sum_n (e_n,A^*B e_n)_H
=\tr(A^*B)$.  
One also has  
$
\|K\| \le \|K\|_2 \le \|K\|_1,
$
see \cite[Thm. V.22]{ReedSimonVol1} or \cite[Section 18]{Conway}. 
For a positive and 
self-adjoint operator, we have
$K \in \CL_1(H)$ if and only 
if $\sqrt{K} \in \CL_2(H)$. To see this note 
 that if $K \ge 0$
and $K=K^*$, then  
$|K|=K=\sqrt{K}\sqrt{K}=\sqrt{K}^*\sqrt{K}$ and hence
$\|K\|_1=\|\sqrt{K}\|_2$. 
We denote with $\CLL(H)$ the closed subset of $\CL_1(H)$ of \emph{self-adjoint, non-negative} trace-class operators. 
 Note that 
 $\CLL(H)$ is closed since 
$\|\cdot\|_1$ is stronger than the operator norm $\|\cdot\|$, and  $\CLL(H)$ is easily seen to be closed with respect to  $\|\cdot\|$.
Below we list some useful facts. 

\begin{enumerate}[label=(K\arabic*),start=2]
\item \label{item:sep_cond} 
The dual of the Banach space $(\CL_1(H),\|\cdot\|_1)$
is isometrically isomorphic to  $(\CL_\infty(H),\|\cdot\|)$, see  \cite[Thm. 19.2]{Conway}. The duality is 
$
\langle A, B \rangle\colon =\tr(A B) $ for $A \in \CL_\infty(H)$ and $B \in \CL_1(H)$.
\item\label{K4} If $A$ and $B$ are Hilbert Schmidt then 
 $\|AB\|_1 \le \|A\|_2\|B\|_2$, 
see \cite[Ex. 28]{ReedSimonVol1} or \cite[Prop. 18.2]{Conway}.
\item  \label{item:Hilbert-Schmidt}   Consider a measure space $(\VV,\mu)$ and  the Hilbert space $H=L^2(\VV,\mu)$.
The space $\CL_2(H)$ of  Hilbert-Schmidt operators
$K: H \to H$ can be identified with 
$L^2(\VV^2,\mu \otimes \mu)$ by $K  \in  \CL_2(H)\longleftrightarrow \CKK \in L^2(\VV^2,\mu \otimes \mu)$
where 
\begin{equation*}\label{L2mapHS}
Kf(y)= \int_{\VV} \CKK(x,y)f(y) \mu(dx),
\qquad \forall \,f \in L^2(\VV,\mu) 
\end{equation*}
and   $\|\CKK\|_{L^2(\VV^2)}
=\| K\|_2$. See, e.g., \cite[Thm. VI.23]{ReedSimonVol1}.
\item \label{Mercer} 
 Le $\VV$ be a compact metric space, $\mu$ a
 Borel measure on $\VV$ and $H=L^2(\VV,\mu)$. Let $\CC^{+,s}$ be the class of continuous, symmetric, Mercer 
(i.e. integrally positive definite) kernels , that is 
\[
\begin{split}
\CC^{+,s}(\VV)= \Big \{ & \CKK \in C^{0}(\VV^2,\RE): \CKK(x,y)=\CKK(y,x) \,\, \forall x,y\in \VV; \, \\
& 
	\inf_{f \in L^2(\VV,\mu)}\int_{\VV^2} \CKK(x,y)f(x)f(y) \mu(dx) \mu(dy) \geq 0.\Big \}. \\
\end{split}
\]
Let $\phi(\CKK)=K$  the corresponding Hilbert-Schmidt operator, that is
$Kf(y)= \int_{\VV} \CKK(x,y)f(x) \mu(dx)$, for all $f\in H$.  
Mercer's Theorem states that 
for every $\CKK \in \CC^{+,s}(\VV)$
one has $\phi(\CKK) \in \CL_{1}^{+,s}$
and 
\[
\tr(\phi(\CKK))=\int_\VV \CKK(x,x) \mu(dx).
\]
See e.g. \cite[Theorem 3.11.7]{Simon_Operatortheory}.  
\item \label{Pow-Stor+Pecc}
  Powers-St{\o}rmer inequalities.
  \begin{itemize} 
  \item The classical Powers-St{\o}rmer inequality (Lemma 4.1 \cite{Powers}) 
  is:  let 
 $K$ and $K'$ be in $\CLL$,   then 
\[
\|\sqrt{K}-\sqrt{K'}\|_2^2 \le \|K-K'\|_1.
\]
\item 
A variant of the  Powers-St{\o}rmer inequality, proved in Prop. 5.12 in \cite{favaro2023quantitative},  is  
\[
\|\sqrt{K}-\sqrt{K'}\|_2 \le 
|\tr(K)-\tr(K')|^{1/2}+\sqrt2
\|K-K'\|_2^{1/4}\min(\tr(\sqrt{K}),\tr(\sqrt{K'}))^{1/2}.
\]
\end{itemize}
\end{enumerate}

\begin{lemma}\label{Lemma_eqvivalent_conv}
Let $(K_n)_n$ and $K$  be  in $\CLL$, with $H$ separable Hilbert. Then
 the following are equivalent:
\begin{enumerate}
 \item $\|\sqrt{K_n}-\sqrt{K}\|_2 \to 0$;
 \item 
 $\|{K_n}-{K}\|_1 \to 0$;
   \item   $\|{K_n}-{K}\|_2 \to 0$
  and $\tr(K_n)=\|{K_n}\|_1 \to \tr(K)= \|{K}\|_1$. 
\end{enumerate}
\end{lemma}

\begin{proof}
Let us first prove that 
(1) yields (2). To prove the claim,  write 
${K_n}-{K}=\frac12(\sqrt{K_n}-\sqrt{K})B_n
+B_n\frac12(\sqrt{K_n}-\sqrt{K})$
with $B_n=(\sqrt{K_n}+\sqrt{K})$. 
Using  \ref{K4}, 
one gets $\|K_n-K\|_1 \le \|\sqrt{K_n}-\sqrt{K}\|_2 
\|\sqrt{K_n}+\sqrt{K}\|_2$ and the claim follows. 
(2) yields (1) by 
 Powers-St{\o}rmer inequality, see \ref{Pow-Stor+Pecc}. Finally, 
by the variant of the Powers-St{\o}rmer inequality recalled in \ref{Pow-Stor+Pecc}
 if $\|K_n- K\|_2 \to 0$ and
$\tr(K_n) \to \tr(K)$  then $\|\sqrt{K_n}-\sqrt{K}\|_2 \to 0$, which is: 
(3) yields (2). 
\end{proof}

\subsection{LLN and LDP for sums of triangular arrays 
of random variables with values in a Banach space}

In what follows, we state and prove a law of large numbers for triangular arrays of random variables with values in a Banach space, this is the core of the proof of Proposition~\ref{prop:LLN_cond}. Under slightly stronger assumptions an LDP holds as well, this is the main result of \cite{Bolthausen1984}, and we state it here for completeness. This is the main building block of the proof of Proposition~\ref{prop:LDP_cond}.\\

Let $E$ be a separable Banach space with norm $\|\cdot\|_E$.
Denote by $E'$ the dual space of $E$ and by $\langle \cdot,\cdot \rangle$
the duality product on $(E,E')$. 
Given a random (measurable) elements $X$ defined on 
a probability space $(\Omega,\CF,\PP)$  with values in $(E,\CBB(E))$, 
let us recall that $X$ is said to be  Bochner integrable 
if $\E[\|X\|_E]<+\infty$. Under this condition 
the expectation of $\E[X]=\int_\Omega X(\omega) \PP(d\omega)$ is a well-defined 
element of $E$, where the integral is the so-called Bochner integral, 
see, e.g., II.2 in \cite{DiestelUhl}.

\begin{proposition}\label{LLNarray} Let
$(X_{n,k})_{n\geq 1; \,k\leq n}$ be a triangular
array of random variables with values in $E$.
	Assume that 
for every $n \ge 1$ the random elements $(X_{n,1},\dots,X_{n,n})$ are i.i.d. and in addition that 
\begin{enumerate}[label=(\alph*)]
\item $X_{n,1}  \stackrel{\CL}{\to}
    X_{\infty} \sim \mu$  when $n \to \infty $(with respect to the strong topology in $E$);
\item\label{assump_LLN_triang} for some $p  >1$ 
\begin{equation}\label{conditionUI}
\sup_n \E[\|X_{n,1}\|_E^{p}]<+\infty.
\end{equation}
\end{enumerate}
Let  $S_n\colon =\frac1{n} \sum_{i=1}^n X_{n,i}$, then 
\[
\lim_{n\to\infty} \E[\|S_n - \E[X_\infty]\|_E]] =0.
\]
\end{proposition}

\begin{proof}
By (a) and Skorokhod theorem, for every $k \geq 1$
there is a sequence $(\tilde X_{n,k})$ 
and a random variable $\tilde X_{\infty,k}$
defined on a suitable probability space 
$(\tilde \Omega_k,  \tilde \CF_k, \tilde  \PP_k)$
such that $\tilde X_{n,k}\stackrel{\CL}{=}X_{n,1}$ for every $n \geq 1$, 
$\tilde X_{\infty,k}\stackrel{\CL}{=}X_{\infty}$ and 
$\tilde X_{n,k} \to \tilde X_{\infty,k}$ with $\tilde \PP_k$ probability one. 
Now, consider $\tilde X_{n,k}$ and $\tilde X_{\infty,k}$ as random elements defined 
on $\tilde \Omega=\prod_{k=1}^{+\infty} \tilde \Omega_k$ and 
$\tilde \CF=\otimes_{k=1}^{+\infty} \CF_k$ with 
the product measure $\tilde \PP=
\otimes_{k=1}^{+\infty} \tilde \PP_k$. In this way   $(\tilde X_{n,k})_{k}$ 
and $(\tilde X_{\infty,k})_{k}$ are independent and 
hence $\tilde S_{n}=\frac{1}{n} \sum_{k=1}^n \tilde X_{n,k} \stackrel{\CL}{=}S_{n}$.  Set also 
$\tilde S_{n}^{\infty}=\frac{1}{n} \sum_{k=1}^n \tilde X_{\infty,k}$. 
Now
\[
\begin{split}
 \E[\|S_n - \E[X_\infty\|_E]]
=\E[\|\tilde S_n - \E[X_\infty\|_E]] 
& \leq \E[\|\tilde S_n-\tilde S_n^{\infty} \|_E] +\E[\|\tilde S_n^{\infty} - \E[X_\infty]\|_E]   \\
& \leq  \E[\|\tilde X_{n,1}-\tilde X_{\infty,1} \|_E]
+\E[\|\tilde S_n^{\infty} - \E[X_\infty]\|_E]. 
\end{split}
\]
Now, by the strong law of large numbers for iid random elements with values in a separable Banach space, since 
$\E[\|\tilde X_{\infty,k}\|_E]=\E[\|X_{\infty}\|_E]<+\infty$,
one has 
$\lim_n \E[\|\tilde S_n^{\infty}-\E[X_\infty]_E\|] \to 0$.
See \cite{Azlarov}. 
On the other hand, $\|\tilde X_{n,1} -\tilde X_{\infty,1}\|_E \to 0$ with $\tilde \PP$ probability one. This, combined with condition \eqref{conditionUI}, which yields uniform integrability of $\|\tilde X_{1,n} -\tilde X_{\infty,1}\|_E$, 
gives  $\lim_{n\to\infty} \E[\|\tilde X_{1,n} -\tilde X_{\infty,1}\|] =0$.   
The thesis follows. 
\end{proof}

\begin{proposition}[Thm. 2 in \cite{Bolthausen1984}]\label{theoBolth}
Let the same assumptions of
Proposition \ref{LLNarray} hold
with \ref{assump_LLN_triang} replaced by the stronger 
\begin{enumerate}[label=(\alph*')]
\setcounter{enumi}{1}
\item\label{ass_b'} for every $t>0$ 
\begin{equation*}\label{Bolt2bis}
\sup_n \E[e^{t \|X_{n,1}\|_{E}}]<+\infty. 
\end{equation*} 
\end{enumerate}
Then, the sequence of laws of $S_n\colon =\frac1{n} \sum_{i=1}^n X_{n,i}$ satisfies a LDP with speed $n$ and good rate function 
\[
h(e|\mu)\colon =\sup_{e' \in E'} \Big \{ \langle e',e \rangle 
- \log\big (\E[e^{ \langle e',X_{\infty} \rangle} ]\big)  
\Big \} \quad X_\infty \sim \mu.
\]
\end{proposition}

Even if it is not explicitly stated in  \cite{Bolthausen1984}, the rate function $h$ appearing in the previous theorem is good. One easy way to see this is showing that $S_n$ is exponentially tight and then apply
\cite[Lemma 1.2.18]{DemboZeitouni}. 
Using \ref{ass_b'} and Proposition \ref{thm:deAcosta}, one gets that there is a convex, compact and well-balanced  set $V$ such that 
$\sup_{n}\E[ e^{q_V (X_{n,1} )  }]  =M<+\infty$. 
Hence, arguing as in the proof of Lemma 
\ref{Lemma:compactness3} 
one gets
\[
 \PP \Big  ( \frac1{ n} \sum_{i=1}^n X_{n,i}  \not \in t_M V   \Big )     \le 
 e^{-t_M n}  \Big ( \E  \big[ e^{q_V (X_{n,1} )  } \big] \Big)^n  \\
  \le  e^{-t_M  n +\log(M)} = e^{-Rn}
\]
for $t_M=R \log(M)$. This shows that 
$S_n$ is exponentially tight.

\subsection{Tail bound for the maximum singular values of
Gaussian matrices} 

Let $\tilde W^{n_1,n_2}$ be a $n_1 \times n_2$ random matrix with $\tilde W_{ij}^{n_1,n_2} \stackrel{iid}{\sim} \CN(0,\la/n_1)$ and 
denote by $\| \tilde W^{n_1,n_2}\|$ its maximum singular values.  
Following the notation in \cite{Ver19}, given a random variable $X$ let 
$\|X\|_{\psi_2}=\inf \{ t >0 : \E[e^{X^2/t^2}] \le 2\}$.
By Ex. 2.5.8 (a) in \cite{Ver19}, $\|\tilde W^{n_1,n_2}_{ij}\|_{\psi_2} \le c_0 \frac{\sqrt{\la}}{\sqrt{n_1}}$ with $c_0$ an absolute constant. Hence, a direct application of
 \cite[Theorem 4.4.5]{Ver19}  gives the following.
\begin{proposition}\label{thmVer19}
Let $\tilde W^{n_1,n_2}$ be a $n_1 \times n_2$ random matrix with $\tilde W_{i,j}^{n_1,n_2} \stackrel{iid}{\sim} \CN(0,\la/n_1)$. Then,
    There is a constant $C=C(\la)$ such that 
    for every $n_1$ and  $n_2$
    \[
    P\Big\{ \|\tilde W^{n_1,n_2}\| > C \Big (\sqrt{\frac{n_2}{n_1}} +1+t\Big) \Big \} \le  2 e^{-n_1 t^2}.
    \]
\end{proposition}

\subsection{Tilted LDP}
We give here a version of Varadhan's Lemma~\cite{Var66} that allows to derive an LDP for a sequence of measures coming from the tilting of a sequence of measures satisfying an LDP. 

\begin{proposition}\label{VaradhanModified}
   Let $(\BX, \mathcal{X})$ be a Polish space with associated Borel $\s$-algebra. Let $\{P_N\}_N$ be a sequence of probability measures on $(\BX, \mathcal{X})$, satisfying an LDP with speed $N$ and good rate function $I$.  
Assume that 
\begin{itemize}
\item[(i)] $\Phi_0$ is a upper bounded and continuous function 
from $\mathbb{X} \to \RE$;
\item[(ii)] $\rho :\mathbb{X} \to [0,+\infty)$ is locally bounded.
\end{itemize}
Then, the sequence of probability measures   $\{P_N^{\circ}\}_N$ defined as follows
\[
P_N^{\circ}(dx)= \frac{e^{-(N\Phi_0(x)+\rho(x))} P_N(dx)}{
 \int_{\mathcal{\mathbb{X}}}
e^{-(N\Phi_0(s) +\rho(s)) } P_N(ds)} ,
\]
satisfies an LDP with speed $N$ and rate function $I(x)-\Phi_0(x)-I_0$  
where $I_0=\inf_x [I(x)-\Phi_0(x)]$.

\end{proposition}

In this form, the above theorem is a slight extension of the result presented in \cite[Theorem III.17]{denHol00}.

\begin{proof}
    First, let us prove that
    \begin{equation}\label{eq:varadhan}
        \lim_{N\to\infty} \frac 1N \log 
 \int_{\mathcal{\mathbb{X}}}
e^{-(N\Phi_0(s) +\rho(s)) } P_N(ds)=-\inf_x [I(x)-\Phi_0(x)].
    \end{equation}
As a lower bound, fix any $x\in \mathcal{\mathbb{X}}$, than for any $B_x$ neighborhood of $x$, we have that
\begin{align*}
    \int_{\mathcal{\mathbb{X}}}
\frac 1N\log e^{-(N\Phi_0(s) +\rho(s)) } P_N(ds)&\geq \frac 1N\log \int_{B_x}
e^{-(N\Phi_0(s) +\rho(s)) } P_N(ds)\\
&\geq -\sup_{y\in B_x}\Phi_0(y)-\delta-\sup_{y\in B_x}\frac1N\rho(y)-\inf_{y\in B_x}I(y)+o(1),
\end{align*}
where we use the continuity of $\Phi_0$ and the fact that $P_N$ satisfies a large deviation principle. Then, taking the limit $N\to\infty$, since $\rho$ is locally bounded, we obtained that 
\[
\liminf_{N\to\infty}  \int_{\mathcal{\mathbb{X}}}
\frac 1N\log e^{-(N\Phi_0(s) +\rho(s)) } P_N(ds)\geq \Phi_0(x)-I(x),
\]
for any $x\in \mathcal{\mathbb{X}}$. The upper bound is a direct consequence of Varadhan's Lemma, since 
\[
\int_{\mathcal{\mathbb{X}}}
\frac 1N\log e^{-(N\Phi_0(s) +\rho(s)) } P_N(ds)\leq \int_{\mathcal{\mathbb{X}}}
\frac 1N\log e^{-N\Phi_0(s) } P_N(ds),
\]
because of the positivity of $\rho$. Once proved \eqref{eq:varadhan}, the thesis follows the same lines of \cite[Theorem III.17]{denHol00}.
\end{proof}

\subsection{Important notation}\label{appendix:notation}
To help the reader, here we collect the most important pieces of notation used throughout the paper.

\begin{center}
\begin{tabular}{c|l}
\multicolumn{2}{c}{\textbf{Parameters}}\\
  $L$   &  depth of the network\\
   $N_0$ and $N_{L+1}=D$ &  input  and output dimensions \\
   $N_1,\cdots,N_{L}$ & width of each hidden layer\\
 $ \UU \subset \RE^{N_0}$ & compact set of possible inputs \\

\multicolumn{2}{c}{\rule{0pt}{4ex}\textbf{Network variables}}\\
 $\vartheta  = \{W^{(\ell)} \}_{\ell=0}^L$ & collection of all trainable weights \\
  $h^{(\ell)}(\bx)$ & pre-activations of  layer $\ell$ from input $\mathbf{x}$, see 
  \eqref{main_recursion}\\
$\CK^{\ell}_{N_{\ell-1}}(\bx,\bx')$ &  conditional covariance function at layer $\ell$ in $(\bx,\bx')$, see \eqref{eq:K_in_C0}\\
$\bK^\ell_{N_{\ell-1}}$ & conditional covariance operator  at layer $\ell$,  see \eqref{conditional_covarianceOp} \\
  $\{\bx_\mu,\by_\mu\}_{\mu=1}^P$ & training set (input and response/label)  \\
 
\multicolumn{2}{c}{\rule{0pt}{4ex}\textbf{Functional spaces, maps and measures}}\\
 $H=L^2(\UU)$ and $ \|\cdot\|_H$ & space of square integrable functions on $\UU$ and corresponding norm \\
$C^0(\UU,\RE^D),C^0(\UU^2,\RE)$ & spaces of continuous functions \\
 $\CC^{+,s}$ &   space of continuous, symmetric, 
positive definite kernels 
on $\UU^2$
\\
$\CL_1$ and $\|\cdot\|_1$ &  space of trace-class operators on $L^2(\UU)$
and  corresponding trace norm
 \\
 $\tr(K)$ & trace of the operator $K \in \CL_1$ \\ 
$\CL_1^{+,s}$ &  space of {non-negative} and {symmetric} trace-class operators
on $L^2(\UU)$ \\
$\CL_\infty$  &  space of bounded linear operators on $L^2(\UU)$,  dual of  $\CL_1$ \\
$\CN_H(\bzero,K)$ & 
Gaussian distribution on $H=L^2(\UU)$  with zero mean 
and covariance $K \in \CL_1^{+,s}$ \\
 $\CKK \mapsto \phi(\CKK)$ & the continuous map  from $\CC^{+,s}$ to $\CL_1^{+,s}$,  defined by 
\eqref{mapKerneltoCov}
\\
 $f \mapsto C_f$ & the continuous map 
 from $H$ to $\CL_1^{+,s}$, defined by \eqref{def_Cf}
\\ 
$\mathcal Q_N(\cdot)$ &  the prior distribution for $(\CK^2_{N_1},\dots, \CK^{L+1}_{N_L})$
\\
$\mathcal Q_N(\cdot | \by_{1:P})$ & the posterior distribution for $(\CK^2_{N_1},\dots, \CK^{L+1}_{N_L})$
given $\by_{1:p}=[\by_1,\dots,\by_P]$\\ 
\end{tabular}
\captionof{table}{Table of important notation} \label{tab:notations}
\end{center}

\section*{Acknowledgments}
LA acknowledges partial financial support by the Italian Ministry of University and Research (MUR) via PRIN 2022 – ConStRAINeD-CUP-2022XRWY7W and  by the European Union-Next
Generation EU, Missione 4-Componente 1-CUP-D53D23018970001 via the project
``Stochastic particle-based anomalous reaction-diffusion
models with
heterogeneous interaction for radiation therapy'' Prot.~P2022TX4FE\_02.
FB is partially supported by the MUR - PRIN project ``Discrete random structures for Bayesian learning and prediction'' no. 2022CLTYP4.
 CH is grateful for the hospitality during a visit at Politecnico di Milano, where this work was initiated. CH was supported by the research grant (VIL69126) from Villum Fonden.

\bibliographystyle{abbrv}
\bibliography{biblio}

 \end{document}